\documentclass[letterpaper,11pt]{article}
\usepackage[letterpaper,left=1in,right=1in,top=1.5in,bottom=1.5in]{geometry}
\usepackage{bbm}
\usepackage{relsize}

\usepackage{amsmath,amssymb,amsthm}

\newcommand{\myauthor}{}

\newcommand*{\TitleFont}{%
      \usefont{\encodingdefault}{\rmdefault}{}{n}%
      \fontsize{15}{25}%
      \selectfont}

\newcommand{\mytitle}{Zeta Spectral Triples\\
\TitleFont{}}

\title{\mytitle}

\author{Alain Connes, Caterina Consani and Henri Moscovici}
\date{}

\usepackage[pdfstartview=FitH,
            pdfauthor={\myauthor},
            colorlinks=true,
            linkcolor=reference,
            citecolor=citation,
            urlcolor=e-mail,
            backref]{hyperref}

\usepackage[
]{cleveref}

\usepackage{amsmath}
\usepackage{amscd}
\usepackage{amsbsy}
\usepackage{amssymb}
\usepackage{verbatim}
\usepackage{eufrak}
\usepackage{eucal}
\usepackage{microtype}
\usepackage{mathrsfs}
\usepackage{amsthm}
\usepackage{stmaryrd}
\usepackage[all,cmtip]{xy}
\usepackage{tikz}
\usetikzlibrary{matrix,arrows}
\usepackage[abbrev,lite,alphabetic]{amsrefs}
\usepackage{dsfont}

\usepackage{fancyhdr}
\usepackage{color}
\definecolor{todo}{rgb}{1,0,0}
\definecolor{conditional}{rgb}{0,1,0}
\definecolor{e-mail}{rgb}{0,.40,.80}
\definecolor{reference}{rgb}{.20,.60,.22}
\definecolor{mrnumber}{rgb}{.80,.40,0} 
\definecolor{citation}{rgb}{0,.40,.80} 

\let\oldmarginpar\marginpar
\renewcommand\marginpar[1]{\-\oldmarginpar[\raggedleft\footnotesize #1]%
{\raggedright\footnotesize #1}}

\usepackage{float}

\newtheorem{theorem}{Theorem}[section]
 \newtheorem{cor}[theorem]{Corollary}
 \newtheorem{lem}[theorem]{Lemma}
 \newtheorem{prop}[theorem]{Proposition}
 
 \theoremstyle{definition}
 \newtheorem{defn}[theorem]{Definition}
 \theoremstyle{remark}

 \numberwithin{equation}{section}

\newcommand{\ie}{{\it i.e.\/}\ }

\newcommand{\cf}{{\it cf.}}

\newcommand{\fourier}{{\F}}

\newcommand{\cD}{\mathcal D}
\newcommand{\cE}{\mathcal{E}}
\newcommand{\cH}{\mathcal{H}}
\newcommand{\cI}{\mathcal I}

\newcommand{\cM}{\mathcal{M}}

\newcommand{\cW}{{\mathcal W}}
\newcommand{\detreg}{\mathop{\mathrm{det}}\nolimits_{\mathrm{reg}}}

\newcommand{\C}{{\mathbb C}}
\newcommand{\F}{{\mathbb F}}

\newcommand{\N}{\mathbb{N}}
\newcommand{\Q}{{\mathbb Q}}
\newcommand{\R}{\mathbb{R}}

\newcommand{\Z}{\mathbb{Z}}

\newcommand{\dln}{{D_{\log}^{(\lambda,N)}}}
\newcommand{\Det}{{\rm Det}}

\newcommand{\id}{{\rm id}}

\newcommand{\Ker}{{\rm Ker}}

\numberwithin{equation}{section}

\hypersetup{linkcolor=black}
\begin{document}
\maketitle







\begin{abstract}
We propose and investigate a strategy toward a proof of the Riemann Hypothesis based on a spectral realization of its non‑trivial zeros. 
Our approach constructs self‑adjoint operators \( \dln \) obtained as rank‑one perturbations of the spectral triple associated with the scaling operator on the interval \([ \lambda^{-1}, \lambda ]\). 
The construction only involves the Euler products over the primes \( p \leq x = \lambda^2 \) and produces self‑adjoint operators whose spectra coincide, with striking numerical accuracy, with the lowest non‑trivial zeros of \( \zeta(\tfrac12 + i s) \), even for small values of \(x\). 
The theoretical foundation rests on the framework introduced in~\cite{VJ} together with the extension in~\cite{CS} of the classical Carathéodory–Fejér theorem for Toeplitz matrices, which guarantees the necessary self‑adjointness. 
Numerical experiments show that the spectra of the operators \( \dln \) converge towards the zeros of \( \zeta(\tfrac12 + i s) \) as the parameters \( N, \lambda \to \infty \). 
A rigorous proof of this convergence would establish the Riemann Hypothesis. 
We further compute the regularized determinants \( \det\nolimits_{\mathrm{reg}}(\dln - z) \) of these operators and discuss the analytic role they play in controlling and potentially proving the above result by showing that, suitably normalized, they converge towards the Riemann $\Xi$ function.
\paragraph{Mathematics Subject Classification 2020.}
\href{http://www.ams.org/mathscinet/msc/msc2020.html?t=11Mxx&btn=Current}{11M06},
\href{http://www.ams.org/mathscinet/msc/msc2020.html?t=11Mxx&btn=Current}{11M55},
\href{http://www.ams.org/mathscinet/msc/msc2020.html?t=22E46&btn=Current}{58B34},
\href{http://www.ams.org/mathscinet/msc/msc2020.html?t=33D45&btn=Current}{33D60},
\href{http://www.ams.org/mathscinet/msc/msc2020.html?t=34B20&btn=Current}{34B20}.

\end{abstract}

\maketitle


\section{Introduction}
This paper is motivated by the spectral realization initiated in \cite{VJ} of the low lying zeros of the Riemann zeta function, in other words as  the infrared part of the spectrum of a selfadjoint operator. Another key ingredient is the generalization proved in \cite{CS} of a fundamental result on Toeplitz matrices which is a   corollary of Carathéodory-Fejer 1911 structure theorem \cite{CF11}. This generalization provides for us a large class of functions whose zeros are  located on the critical line $\Re(s) = 1/2$ due to the selfadjointness of relevant matrices and the Hurwitz theorem on zeros of  uniform limits of holomorphic functions. With these tools at hand one devises a process, perfectly in line with Riemann's paper \cite{Riemann},  which associates to the restricted  Euler product involving only the primes $p\leq x$ a function whose zeros are on the critical line. \newline
The great surprise then, is that the zeros of this function give high-precision approximations to the first non-trivial zeros of the Riemann zeta function using remarkably few terms of the Euler product. For instance using only the  primes $\leq 13$ one obtains for the first 50 zeros an extraordinary accuracy, with errors ranging from $2.5 \times 10^{-55}$ for the first zero to approximately $10^{-3}$ for the fiftieth.  The probability of achieving such precise approximations by chance is approximately $10^{-1235}$, effectively ruling out coincidence and suggesting a deep structural relationship between the restricted Euler products and the location of the zeros.
 The method we use is general as well as the proof  that all the approximating values lie exactly on the critical line.\newline  
In fact we construct spectral triples associated to rank one perturbations $\dln$ of the scaling operator $D_{\log}^{(\lambda)}$ on the interval $[\lambda^{-1},\lambda]$ with periodic boundary conditions. This construction is based on the restriction $QW_\lambda^N$ of the Weil quadratic form to the linear space $E_N=E_N(\lambda)$ of test functions spanned by the $2N+1$ eigenfunctions of $D_{\log}^{(\lambda)}$ associated to the $2N+1$ eigenvalues of smallest absolute value (\ie $\leq N \pi/\log \lambda$),  extended by $0$ outside the interval $[\lambda^{-1},\lambda]$. We need to verify that the smallest eigenvalue $\epsilon_N$ of $QW_\lambda^N$ is simple and that the corresponding eigenfunction is "even" \ie invariant under the symmetry $u\mapsto u^{-1}$. We let $\delta_N\in E_N$ be the vector representing  the Dirichlet kernel, as an approximation to the evaluation on the boundary of the interval $[\lambda^{-1},\lambda]$. Our main result is the following:

\begin{theorem}\label{finmainintro} Let $\epsilon_N$ be the smallest eigenvalue of $QW_\lambda^N$  assumed simple and  $\xi$ the corresponding eigenvector assumed even, normalized by $\delta_N(\xi)=1$.\newline
 $(i)$~The operator $\dln=D_{\log}^{(\lambda)}-\vert D_{\log}^{(\lambda)}\xi\rangle\langle \delta_N\vert$ is selfadjoint in the direct sum $E'_N\oplus E_N^\perp$ where  on the subspace $E'_N=E_N/\C\xi$ the inner product is given 	by the restriction of the quadratic form $QW_\lambda^N-\epsilon_N \langle \mid \rangle$.\newline 
 $(ii)$~The regularized determinant of $\dln$ is given by $\detreg(\dln - z)=-i\,\lambda^{-iz}\widehat\xi(z)$ where $\widehat\xi$ is the Fourier transform of $\xi$ for the duality $\langle\R_+^*\mid \R\rangle$.\newline
 $(iii)$~The Fourier transform $\widehat\xi(z)$ is an entire function, all its zeros are on the real line and coincide with the spectrum of $\dln$.
 \end{theorem}
 In Section \ref{sectnumerical} we show the striking numerical evidence for the convergence of the eigenvalues of the selfadjoint operators $\dln$ towards the zeros of the Riemann zeta function $\zeta\left(\frac{1}{2}+i s\right)$ as the parameters $N$ and $\lambda$ tend to infinity.\newline
 In Section \ref{outlook}, we explain the natural strategy to justify the above numerical convergence. It consists in taking the first steps in trying to show the convergence of the regularized determinants $\detreg(\dln - s)$ towards the Riemann $\Xi$ function.\newline
 Finally Section \ref{missing} describes the missing 
 steps in the above strategy and the perspective that it 
  opens, based on \cite{Co-zeta}, on  the connections   
  between the world of the Weil quadratic form and that  
  of information theory, as developed through the theory   
   of prolate wave functions by D. Slepian and his collaborators \cite{Slepian}.
 
\section{Preliminaries}
\subsection{The Banach algebra $L^1(\R,dx)$}
In this section we shall explain elementary computations preparing the ground for the introduction of the Weil quadratic form $QW_\lambda$ and of the matrix encoding this quadratic form in a natural  basis. 
\begin{defn} We let $L^1(\R,dx)$ be the involutive complex Banach algebra of complex valued integrable functions on $\R$ with product and involution given by 
\begin{equation}\label{convo}
	f*g(y):=\int f(x)g(y-x)dx
\end{equation}
\begin{equation}\label{pinvolution}
	f^*(y):=\overline{f(-y)}
\end{equation}	
\end{defn}
We shall use the subalgebra of compactly
supported functions. Both operations yield functions with compact support. 
We consider the inclusion $L^2([0, L])\subset L^1(\R,dx)$ obtained by extending functions by $0$ outside the interval $[0,L]$.
In particular, for $f, g \in L^2[0, L]$ the support of $f^* \ast g$ is contained in $[-L, L]$.\newline 
We shall use inner products $\langle f,g\rangle$ which are antilinear in the first variable and linear in the second. The standard inner product $\langle f,g\rangle_2$  is given by 
\begin{equation}\label{inner}
	\langle f,g\rangle_2=\int \overline{f(x)}g(x)dx=(f^**g)(0)
\end{equation}
In preparation for more involved inner products, we shall compute expressions of the form 
\begin{equation}\label{qfg}
	q(f,g)(y):=(f^**g)(y)+(f^**g)(-y)
\end{equation}
which is an even function of $y$ depending antilinearly on $f$ and  linearly on $g$. One has 
$$
(f^**g)(-y)=\overline{(f^**g)^*(y)}=\overline{(g^**f)(y)}, 
$$
hence 
\begin{equation}\label{qfg1} q(f,g)(y)=(f^**g)(y)+\overline{(g^**f)(y)}.
\end{equation}
Let $a\in \R$ and for any $f\in L^1(\R,dx)$, let $f_a(x):=f(x-a)$ denote the translated of $f$.
\begin{lem}\label{qtrans} For any $a\in \R$, $f,g\in L^1(\R,dx)$, one has $(f_a)^**g_a=f^**g$ and $q(f_a,g_a)=q(f,g)$.	
\end{lem}
\proof One has 
$$
(f_a)^**g_a(y)=\int (f_a)^*(y-x)g_a(x)dx=$$ $$=\int \overline{(f_a)(x-y)}g_a(x)dx=
\int \overline{f(x-a-y)}g(x-a)dx
$$
which is independent of $a$. The second equality follows from \eqref{qfg}. \endproof 
\subsection{The basis $\{U_n\}_{n \in \Z}$}
 
A natural orthonormal basis for $L^2([0, L])$ consists of the functions $\{U_n\}_{n \in \Z}$  defined by
\begin{align} \label{expbasis}
U_n(x):=L^{-\frac 12} \exp(2\pi i n x/L) , \qquad \forall x	\in [0,L] .
\end{align}
Applying the above operations to the basis functions one obtains, for $n\neq m$, $y\in [0,L]$, \begin{small}
$$
(U_m^*\ast U_n)(y)= \int \overline{ U_m(x-y)}U_n(x)dx=\frac 1L\int_y^L\exp\big(2\pi im(y-x)/L+2\pi inx/L\big)dx $$
$$=\frac{\exp(2\pi im y/L)}{L}\int_y^L\exp(2\pi i(n-m)x/L)dx=\frac{\exp(2\pi im y/L)}{2\pi i(n-m)}
\big(\exp(2\pi i(n-m)x/L) \big)\big\vert_y^L$$                              \label{sim-mn}
$$=\frac{\exp(2\pi im y/L)}{2\pi i(n-m)}\left(1-\big(\exp(2\pi i(n-m)y/L) \big)\right)
=\frac{\exp(2\pi im y/L)-\exp(2\pi in y/L)}{2\pi i(n-m)} 
$$
\end{small}
so that for $n\neq m$, $y\in [0,L]$, 
\begin{equation}\label{Umn}
	(U_m^*\ast U_n)(y)=\frac{\exp(2\pi im y/L)-\exp(2\pi in y/L)}{2\pi i(n-m)}
\end{equation}

The result of \eqref{Umn} being symmetric in $n,m$, it implies using \eqref{qfg1}, that for $y\in [0,L]$,
\begin{align}
(U_m^**U_n)(y)+(U_m^**U_n)(-y)&=2\,\Re\left(\frac{e^{2\pi im y/L}-e^{2\pi in y/L}}{2\pi i(n-m)}\right)\nonumber\\ &=	\frac{\sin(2\pi m y/L)-\sin(2\pi n y/L)}{\pi (n-m)} \label{symrealmn}
\end{align}
so that the even function $q(U_m,U_n)$ is given by 
\begin{equation}\label{symrealmn1}
q(U_m,U_n)(y)=	\frac{\sin(2\pi m \vert y\vert /L)-\sin(2\pi n \vert y\vert/L)}{\pi (n-m)}, \  \  \forall y\in [-L,L].
\end{equation}
For $m=n$ one simply has, for $y\in [0,L]$,
$$
(U_n^**U_n)(y)=\frac 1L\int_y^L\exp(2\pi in(y-x)/L+2\pi inx/L)dx=(1-y/L)\exp(2\pi in y/L) ,
$$
whence
$$
(U_n^**U_n)(y)+(U_n^**U_n)(-y)=2\,\Re\left((1-y/L)\exp(2\pi in y/L)\right)=$$
$$=2(1-y/L) \cos(2\pi n y/L).$$
We thus get 
\begin{equation}\label{symrealn1}
q(U_n,U_n)(y)=2(1-\vert y\vert /L) \cos(2\pi n y/L) , \  \  \forall y\in [-L,L].
\end{equation}
\begin{lem}\label{polarize0}
The explicit expression of the even function $q(U_n,U_m)(y)$ for $y\in [0,L]$ 
is :
$$
 \begin{array}{ccc}
 	 	m\neq n &\vline & \frac{\sin \left(\frac{2 \pi  m y}{L}\right)-\sin \left(\frac{2 \pi  n y}{L}\right)}{\pi(n-m)}\\
 	m= n &\vline & \ \frac{2(L-y) \cos \left(\frac{2 \pi  n y}{L}\right)}{L}
 \end{array}
$$
\end{lem}
\section{The Weil quadratic form $QW$} 
The sesquilinear form $QW$ derives from Weil's formulation of the explicit formula in prime numbers 
theory \cite{Weil}, which we recall below. Let us denote by $\cW(\R^*_+)$ the 
  {\em Weil  class}  of complex valued functions $f$ on $\R^*_+$, by which we mean the functions
having continuous derivative, except at finitely many points where both $f(x)$ and $f'(x)$ may
have at most a discontinuity of the first kind; at such a point the value of $f(x)$ and $f'(x)$ is defined as the average of the right and left limits. In addition, the functions $f \in W(\R^*_+)$ are assumed to satisfy the estimate
$$
f(x)=O(x^\delta), \ \text{for} \ x\to 0+, \ \ f(x)=O(x^{-1-\delta}), \ \text{for} \ x\to +\infty ,
$$  
 for some  $\delta>0$. These functions admit a Mellin transform, denoted 
\begin{equation}\label{mellin}
 \tilde f(s):=\int_0^\infty f(x)x^{s-1}dx .
 \end{equation}
With the additional notation \  $f^\sharp(x):=x^{-1}f(x^{-1})$, the Weil's explicit formula 
takes the form (\cf  \cite{BombieriRH})
 \begin{equation}\label{bombieriexplicit}
 \sum_{\rho}\tilde f(\rho)=\int_0^\infty f(x)dx+\int_0^\infty f^\sharp(x)dx-\sum_v {\mathcal W}_v(f),
 \end{equation}
 where $\rho$ runs over all complex zeros $\rho$ of the Riemann zeta function,  $v$
 runs over all rational places of $\Q$,  the non-archimedean distributions $\mathcal W_p$ are defined as 
 \begin{equation}\label{bombieriexplicit1}
 {\mathcal W}_p(f):=(\log p)\sum_{m=1}^\infty\left(f(p^m)+f^\sharp(p^m)\right),
 \end{equation}
and the archimedean distribution is given by 
 \begin{equation}\label{bombieriexplicit2}
 {\mathcal W}_\R(f):=(\log 4\pi +\gamma)f(1)+\int_{1}^\infty\left(f(x)+f^\sharp(x)-\frac 2x f(1)\right)\frac{dx}{x-x^{-1}} .
 \end{equation}
 It should be noted that the sum in the left hand side of \eqref{bombieriexplicit},
 whose general term is oscillatory, is only conditionally convergent. This issue of lack of absolute convergence is an essential feature of Riemann's formula for the function $\pi(x)$ which is the number of primes less than $x$. In our situation this issue does not appear since we shall only apply the explicit formula to functions which are the convolution (for the group $\R_+^*$) of two square integrable functions with compact support, thus ensuring the absolute convergence of the sum over the zeros.
 
 An equivalent formulation, known as Guinand-Weil formula, uses the Fourier transform 
 \begin{equation} \label{mhat}
 {\widehat F}(s):=\int_{\R^*_+} F(u)u^{-is}d^*u , \qquad  d^*u=\frac{du}{u}
\end{equation}
 in place of the Mellin transform \eqref{mellin}. The passage from one to the other
 is be obtained  by implementing the  automorphism  
 \begin{equation}\label{MF1}
f\mapsto \Delta^{1/2}f=F, \quad \ie \quad F(x)= x^{1/2}f(x) ,
\end{equation}
which respects the convolution product and satisfies the equalities
$$
(\Delta^{1/2}f^\sharp)(x)=x^{1/2}f^\sharp(x)=x^{-1/2}f(x^{-1})=(\Delta^{1/2}f)(x^{-1}).
$$
For a rational place $v$,  denoting
$W_v(F):={\mathcal W}_v(\Delta^{-1/2}F),$  the above distributions $\mathcal W_p$ take the following  form:
\begin{equation}\label{bombieriexplicit1bis}
 W_p(F)=(\log p)\sum_{m=1}^\infty p^{-m/2}\left(F(p^m)+F(p^{-m})\right) ,
 \end{equation}
 while the archimedean distribution $\mathcal W_\R$ becomes
$$
 W_\R(F):=(\log 4\pi +\gamma)F(1)+\int_{1}^\infty\left(F(x)+F(x^{-1})-2x^{-1/2} F(1)\right)\frac{x^{1/2}}{x-x^{-1}}d^*x ,  
 $$  
 where  $ d^*x=dx/x$. The latter can also be expressed as $W_\R = - W_\infty$, where
 \begin{equation}\label{thetaprime}
 W_\infty(F) = \int_\R {\widehat F}(t)\frac{2\partial_t\theta(t)}{2 \pi}dt .
 \end{equation}
 and 
  \begin{equation}\label{riesie}
\theta(t) = - \frac{t}{2} \log \pi + \Im \log \Gamma \left(
\frac{1}{4} + i \frac{t}{2} \right) 
\end{equation}
is  the angular Riemann-Siegel function, with $\log \Gamma(s)$ for $\Re(s)>0$ denoting the 
branch of the logarithm which is real for $s$ real. \newline

By polarization, the Weil form gives the  sesquilinear expression
\begin{align}\label{bombtest}
QW(f,g)=\Psi(f^**g), \ \  
 \Psi(F):= W_{0,2}(F) - W_\R(F)-\sum_{p} W_p(F). 
 \end{align}
The components $W_\R$ and $W_p$ are as above, and the functional $W_{0,2}$ is  \begin{align}\label{bombtest-0}
 W_{0,2}(F) = {\widehat F}(i/2)+{\widehat F}(-i/2) .
 \end{align}

There is a rather subtle invariance property of the Weil sesquilinear form, namely
its symmetry under the inversion $\iota(u)=u^{-1}, \, u \in\R^*$, which 
will play a significant role in its explicit description. 
\begin{lem}\label{wsharp} The Weil functional $\Psi$ fulfills 
\begin{equation}\label{weilQexpw}
 \Psi(h)=\Psi^\#(h)+\Psi^\#(h\circ\iota)=\Psi^\#(h+h\circ\iota) ,
\end{equation}
 where $\Psi^\sharp$ is the distribution on $[1,\infty)$,
\begin{align}\label{psipsisharp}
\Psi^\sharp:=W_{0,2}^\# - W_\R^\#-\sum W_p^\#
 \end{align}
	with the components given by
\begin{align}\label{bombtest02}
W_{0,2}^\#(F)=\int_{1}^{\infty} F(x)(x^{1/2}+x^{-1/2})d^*x ,\\ \label{bombtestR} 
W_\R^\#(F)=\frac 12(\log 4\pi +\gamma)F(1)+\int_{1}^\infty\frac{x^{1/2}F(x)-F(1)}{x-x^{-1}}d^*x ,
\\ \label{bombtestp}
W_p^\#(F)=(\log p)\sum_{m=1}^\infty p^{-m/2}F(p^m).
 \end{align}
\end{lem}
\proof This follows from the construction of $\Psi$. Note the factor $\frac 12$ in \eqref{bombtestR}. \endproof 
\begin{prop}\label{toadd} Let $\lambda>1$, $L=2\log \lambda$.\newline
$(i)$~The following map is an isometry $\kappa:L^2([0,L],dx)\to L^2([\lambda^{-1}, \lambda],d^*u)$,
\begin{equation}\label{toa}
\kappa(f)(u)=f(\log(\lambda u))
\end{equation}
which induces an isomorphism  $C^\infty([0,L])\to C^\infty([\lambda^{-1}, \lambda])$.\newline
	$(ii)$~Let $f,g\in C^\infty([0,L])$. One has 
\begin{equation}\label{weilQexp}
QW(\kappa(f),\kappa(g))=\Psi^\sharp(F), \  \  F(u)=q(f,g)(\log u).
\end{equation}
\end{prop}
\proof $(i)$~Follows since the map $u\mapsto \log(\lambda u)$ is a diffeomorphism from $[\lambda^{-1}, \lambda]$ to $[0,L]$ transforming the measure $d^*u$ into the measure $dx$.\newline 
$(ii)$~One has, by Lemma \ref{wsharp}, 
$$
QW(\kappa(f),\kappa(g))=\Psi^\sharp(h+h\circ \iota), \ \ h=\kappa(f)^**\kappa(g).
$$
Let us show that $h+h\circ \iota=F$. By Lemma \ref{qtrans}, one has, with $a=-\frac L2$ 
$$
q(f,g)=q(f_a,g_a), \ \ f_a(x):=f(x+\frac L2), \ g_a(x):=g(x+\frac L2)
$$
so that $f_a$ and $g_a$ have support in $[-\frac L2,\frac L2]$ and one has $\kappa(f)=f_a\circ \log$, $\kappa(g)=g_a\circ \log$, where  the $\log$ is the  isomorphism of locally compact groups $\log:\R_+^*\to \R$. This  induces an isomorphism of involutive convolution algebras 
$$
\circ \log :L^1(\R,dx)\to L^1(\R_+^*,d^*x).
$$
Thus one gets 
$$
h=\kappa(f)^**\kappa(g)=(f_a\circ \log)^**(g_a\circ \log)=(f_a^**g_a)\circ \log 
$$
which using \eqref{qfg} gives 
$$
h+h\circ \iota=q(f_a,g_a)\circ \log=q(f,g)\circ \log =F.
$$
and hence the required equality.
\endproof 
\subsection{The quadratic form $QW_\lambda$} 
Let $\lambda>1$. We denote by $QW_\lambda$ the restriction of the quadratic form $QW$ to
$L^2([\lambda^{-1}, \lambda],d^*u)$, where $d^*u=\frac{du}{u}$. One has
\begin{equation}\label{quadratsemi}
 QW_\lambda(f,f)=\int_\R 
\vert{\widehat f}(t)\vert^2\frac{2\partial_t\theta(t)}{2 \pi}dt + 2 \Re\left({\widehat f}(\frac i2)\bar {\widehat f}(-\frac i2)\right)-\sum_{1<n\leq \lambda^2} \Lambda(n)\langle f\mid T(n)f\rangle ;
\end{equation}
here $\Lambda(n)$ is the von Mangoldt function, and  $T(n)$ is the bounded self-adjoint operator in $L^2([\lambda^{-1},\lambda],d^*u)$ defined by
\begin{equation}\label{quadratsemi1}
 \langle f\mid T(n)g\rangle=n^{-1/2}\left((f^**g)(n)+(f^**g)(n^{-1})\right ).
\end{equation}
 
\begin{prop} (\cite[\S 2]{VJ}) \label{Hilbert} 
The quadratic form $QW_\lambda$ is lower bounded and lower semi-continuous. 
\end{prop}

Recall that a lower bounded, lower semi-continuous (lsc) quadratic form $Q$ on a Hilbert space $\cH$ is a lower semi-continuous map  $Q:\cH\to (-\infty,+\infty]$, 
\ie such that $Q(\xi)\leq \liminf Q(\xi_n)$ when $\xi_n\to \xi$, which  fulfills
 $Q(\lambda \xi)=\vert \lambda\vert^2 Q(\xi)$ for all $\lambda \in \C$, satisfies
 the parallelogram law 
$$
Q(\xi + \eta) +Q(\xi-\eta)=2 Q(\xi) +2 Q(\eta)
$$
and also an inequality of the form $Q(\xi)\geq -c \Vert \xi\Vert^2$ for all $\xi \in \cH$  reflecting the lower bound of $q$.  The associated sesquilinear form (antilinear in the first variable) is given on the domain ${\rm Dom}(Q):=\{\xi \in \cH\mid Q(\xi)<\infty\}$  by 
  $$
  Q(\xi,\eta):=\frac 14\left( Q(\xi + \eta) -Q(\xi-\eta)+iQ(i\xi + \eta) -iQ(i\xi-\eta)\right).
  $$

 Let $V_n: [\lambda^{-1},\lambda] \to \C$ be the function $\kappa(U_n)$, \ie  
  \begin{equation} \label{vN}
 V_n(u):=U_n(\log (\lambda u)), \  \ \forall u\in [\lambda^{-1},\lambda]
 \end{equation}
  and let $E\subset L^2([\lambda^{-1},\lambda],d^*u)$ be the linear
subspace  generated by the $V_n$ for $n\in \Z$.

  \begin{prop}  (\cite[Prop. 2.3]{VJ}) \label{Hilbert1}  
 The space $E$ is a core for the 
  quadratic form $QW_\lambda:L^2([\lambda^{-1},\lambda],d^*u)\to (-\infty,+\infty]$, 
  which satisfies, for any $\ f\in L^2([\lambda^{-1},\lambda],d^*u)$,
 \begin{equation}\label{quadratsemi7}
 QW_\lambda(f,f)=\liminf_{g_n\to f} QW_\lambda(g_n, g_n), \ g_n\in E .
\end{equation}
In particular, the lower bound of $QW_\lambda$ is the limit, when $N\to \infty$, of the smallest eigenvalue of the restriction of $QW_\lambda$  to the linear span $E_N$ of the functions $V_k$ 
with $\vert k\vert \leq N$.	
\end{prop}

 \subsection{Discrete spectrum of the semilocal Weil quadratic form \texorpdfstring{$Q W_\lambda$}{QWlambda}}
At this juncture we appeal to a basic result from the general theory of quadratic forms. Adopting the notation in \cite[Ch. 10]{sm}, to a quadratic form $\mathfrak{t}$ one associates 
a mapping $\mathfrak{t}^{\prime}: \mathcal{H} \rightarrow \mathbb{R} \cup\{+\infty\}$ by setting $\mathfrak{t}^{\prime}[x]=\mathfrak{t}[x]$ if $x$ is in $\mathcal{D}(\mathfrak{t})$ and $\mathfrak{t}^{\prime}[x]=+\infty$ if $x$ is not in $\mathcal{D}(\mathfrak{t})$. According to 
\cite[Proposition 10.1]{sm} the following four conditions are equivalent:
	\begin{enumerate}
		\item $\mathfrak{t}$ is closed.
		\item If $\left(x_n\right)_{n \in \mathbb{N}}$ is a sequence from $\mathcal{D}(\mathfrak{t})$ such that $\lim _{n \rightarrow \infty} x_n=x$ in $\mathcal{H}$ for some $x \in \mathcal{H}$ and $\lim _{n, k \rightarrow \infty} \mathfrak{t}\left[x_n-x_k\right]=0$, then $x \in \mathcal{D}(\mathfrak{t})$ and $\lim _{n \rightarrow \infty} \mathfrak{t}\left[x_n-x\right]=0$.
		\item $\mathfrak{t}^{\prime}$ is a lower semicontinuous function on $\mathcal{H}$.
		\item If $\left(x_n\right)_{n \in \mathbb{N}}$ is a sequence from $\mathcal{D}(\mathfrak{t})$ such that $\lim _{n \rightarrow \infty} x_n=x$ in $\mathcal{H}$ for some $x \in \mathcal{H}$ and the set $\left\{\mathfrak{t}\left[x_n\right]: n \in \mathbb{N}\right\}$ is bounded, then we have $x \in \mathcal{D}(\mathfrak{t})$ and $\mathfrak{t}[x] \leq \liminf _{n \rightarrow \infty} \mathfrak{t}\left[x_n\right]$.
	\end{enumerate}
 
By Proposition \ref{Hilbert} the condition (3) is  fulfilled by the semilocal Weil quadratic form 
$Q W_\lambda$, and therefore Theorem 10.7 in \cite{sm}, reproduced below, applies.

{\bf Representation theorem for semibounded forms --}
{\em If $\mathfrak{t}$ is a densely defined lower semibounded closed form on $\mathcal{H}$, then the operator $A_{\mathfrak{t}}$ is self-adjoint, and $\mathfrak{t}$ is equal to the form $\mathfrak{t}_{\left(A_{\mathfrak{t}}\right)}$ associated with $A_{\mathfrak{t}}$}.

\
 
Thus, for each $\lambda>1$, there is a canonical lower bounded unbounded selfadjoint  operator $A_\lambda$ in the Hilbert space  $L^2\left(\left[\lambda^{-1}, \lambda\right], d^* u\right)$ such that 
\begin{equation}
	Q W_\lambda(f, f)=\langle A_\lambda f\mid f\rangle .
\end{equation}

By construction the unbounded selfadjoint  operator $A_\lambda$ is lower bounded. The issue is to show that it has discrete spectrum. We use the following from \cite{sm}, Proposition 10.6,
\begin{prop}\label{compres}
	Suppose that $A\geq m_A$ is a lower semibounded self-adjoint operator and $m<m_A$. Then the following assertions are equivalent:
	\begin{enumerate}
	\item	The embedding map $\mathcal{I}_{\mathfrak{t}_A}:\left(\mathcal{D}[A],\|\cdot\|_{\mathfrak{t}_A}\right) \rightarrow(\mathcal{H},\|\cdot\|)$ is compact.
	\item The resolvent $R_\lambda(A)$ is compact for one, hence for all, $\lambda \in \rho(A)$.
	\item $(A-m I)^{-1 / 2}$ is compact.
	\item A has a purely discrete spectrum.
	\end{enumerate}
\end{prop}
\begin{theorem}\label{thmsmallest} The selfadjoint operator $A_\lambda$ has discrete lower bounded spectrum.	
\end{theorem}
\proof By the proof of the lower boundedness in \cite{VJ}, the contribution of the non-archimedean primes to the operator $A_\lambda$ is bounded as well as the contribution of the evaluation of the Fourier transform at the poles. Thus it is enough to deal, for any $\lambda>1$  with the contribution of the archimedean place to $A_\lambda$ in the Hilbert space  $L^2\left(\left[\lambda^{-1}, \lambda\right], d^* u\right)$. It is given, after Fourier transform, by the multiplication by  
\begin{equation}
 \partial_t \theta(t)=\frac{1}{2}(\log (|t|)-\log (2)-\log (\pi))-\frac{1}{48 t^2}+O\left(t^{-4}\right)
\end{equation}
 whose asymptotic expansion allows one to use instead the operator $\mathscr{L}$ of multiplication (in Fourier) by the function which is $1$ for $\vert t\vert \leq e$ and $\log (|t|)$ otherwise. Note that the factor $\frac{1}{2 \pi}$ is taken care of since the unitary Fourier transform has a factor $\frac{1}{\sqrt {2 \pi}}$.\newline
By Proposition \ref{compres}, it is enough to show that the embedding map $\mathcal{I}:\left(\mathcal{D}[\mathscr{L}],\|\cdot\|_{\mathfrak{t}_L}\right) \rightarrow(\mathcal{H},\|\cdot\|)$ is compact. The map $\cI$ is of norm $\leq 1$ by construction and it is enough to show that  the image of the unit ball is precompact in the following sense:\newline
 Let $E$ be a metric space. If any of the following three properties is satisfied, then all three are satisfied, and $E$ is said to be precompact:\newline
1. For every $\varepsilon>0$, $E$ can be covered by a finite number of balls of radius $\varepsilon$;\newline
2. For every $\varepsilon>0$, $E$ can be covered by a finite number of subsets with diameter less than \(\varepsilon\);\newline
3. Every sequence in  $E$ has a Cauchy subsequence.\newline 
 We shall  show that for any positive increasing function $\rho:[0,\infty)\to [1,\infty)$ such that $\rho(u)\to \infty$ when $u\to \infty$ the embedding $I_\rho$ of the Hilbert space $\cD_\rho$ in $L^2\left(\left[\lambda^{-1}, \lambda\right], d^* u\right)$ is compact, where the norm square in $ \cD_\rho$ is given by 
\begin{equation}\label{mfour}
 \Vert f\Vert_\rho^2:= \int|\hat{f}(t)|^2 \frac{\rho(t)}{2 \pi} d t
 \end{equation}
 In fact one can use the logarithm-exponential isomorphism and replace the Hilbert space $L^2\left(\left[\lambda^{-1}, \lambda\right], d^* u\right)$  by $\cH=L^2\left(\left[-L, L\right], d x\right)$ and use the ordinary Fourier transform in \eqref{mfour}. It is enough to show that the image  $I_\rho(B)$ of the unit ball $B$ of $ \cD_\rho$ is precompact in $\cH$. Let then $\epsilon >0$ and let us show that one can cover $I_\rho(B)$ by finitely many balls of radius $\epsilon$ for the norm of $\cH$. Since $\rho(u)\to \infty$ when $u\to \infty$ there exists $T<\infty$ such that 
\begin{equation}\label{four1}
  \Vert f\Vert_\rho^2\leq 1\Rightarrow \int_{\vert t\vert \geq T}|\hat{f}(t)|^2 \frac{dt}{2 \pi}\leq (\epsilon/4)^2 
\end{equation}
 Next, the operator $\widehat{P_T}P_L$ is a compact operator in $L^2(\R)$ and hence the image $\widehat{P_T}P_L(C)$ of the unit ball $C$ of $L^2(\R)$ is precompact in $L^2(\R)$. By construction $I_\rho(B)\subset \cH=L^2\left(\left[-L, L\right], d x\right)=P_L L^2(\R)$, thus the map 
 $$
 \widehat{P_T}I_\rho:  \cD_\rho\to L^2(\R)
 $$ 
 is compact. 
  Thus there exists a finite set of functions 
$\{f_j\}_{j\in J}\subset B\subset \cD_\rho$,  such that 
$$
\forall f\in B, \exists j \mid \Vert \widehat{P_T} (f-f_j)\Vert^2 \leq (\epsilon/2)^2 .
$$
Thus one has 
$$
\forall f \in B, \exists j \in J \quad \text{such that} \quad
  \int_{\vert t\vert \leq T}|\hat{f}(t)-\hat f_j(t)|^2 \frac{dt}{2 \pi}\leq (\epsilon/2)^2 
$$
By \eqref{four1} it follows that for any $f\in B$ there exists $j$ such that 
\begin{align*}
\Vert f-f_j\Vert^2_\cH&=\int|\hat{f}(t)-\hat f_j(t)|^2 \frac{dt}{2 \pi}=\\
=\int_{\vert t\vert \leq T}|\hat{f}(t)&-\hat f_j(t)|^2 \frac{dt}{2 \pi}+\int_{\vert t\vert \geq T}|\hat{f}(t)-\hat f_j(t)|^2 \frac{dt}{2 \pi}  \leq (\epsilon/2)^2+(\epsilon/2)^2<\epsilon^2 ,
\end{align*}
where we used  \eqref{four1} and the triangle inequality to get 
$$
\left(\int_{\vert t\vert \geq T}|(\hat{f}(t)-\hat f_j)(t)|^2 \frac{dt}{2 \pi}\right)^{1/2}\leq 
\epsilon/4+\epsilon/4=\epsilon/2 .
$$
\endproof 
\begin{cor}\label{lowest} Let $\lambda>1$. There exists an element $\phi\in L^2\left(\left[\lambda^{-1}, \lambda\right], d^* u\right)$ such that $A_\lambda(\phi)=\mu_\lambda\, \phi$ where $\mu_\lambda$ is the largest lower bound of the spectrum of $A_\lambda$.	
\end{cor}
Note that we cannot assert that $\mu_\lambda\geq 0$.  One has however
\begin{equation}\label{monotone}
 \lambda > \lambda' \Rightarrow  \mu_\lambda\leq  \mu_{\lambda'}
\end{equation}
Indeed first note that the test functions which are piecewise smooth form a core for the Weil quadratic form $Q W_\lambda$ since the smooth ones already do by \cite{VJ} and moreover the Fourier transform of piecewise smooth functions $f$ with compact support are  $O(\vert s\vert^{-1})$ since the derivative $f'$ is a bounded measure, so the piecewise smooth functions are in the domain of $Q W_\lambda$. One then uses the equivalence (with $f$ piecewise smooth)
$$
\nu \leq \mu_\lambda\iff Q W_\lambda(f,f) \geq \nu \Vert f\Vert^2, \  \ \forall f \mid \ {\rm support}(f) \subset \left[\lambda^{-1}, \lambda\right]
$$
\begin{cor}\label{strange} If the limit when $\lambda\to \infty$ of the decreasing function $\mu_\lambda$ is equal to $0$  then RH  holds.
\end{cor}

\section{The matrix of $QW_\lambda$ in the basis $V_n$} 
We now compute the matrix of the sesquilinear form $QW_\lambda(f,g)= \Psi(f^*\ast g)$ in the basis 
$\{V_n\}_{n \in \Z}$. By Proposition \ref{toadd} and the equality $V_n=\kappa(U_n)$ we get  
\begin{equation}\label{weilformbase}
QW_\lambda(V_n,V_m) =
 \Psi^\#(F), \quad \text{where} \quad
 F(x)=q(U_n,U_m)(\log x)
\end{equation}
and $\Psi^\sharp$ was defined in \eqref{psipsisharp}. \newline
Next, we proceed to describe the contribution to the matrix $QW_\lambda(V_n,V_m)$ of each term in \eqref{psipsisharp}. 

\subsection{The matrix $ W_{0,2}(V_n,V_m)$}
The following lemma  shows that  the terms $ W_{0,2}(V_n,V_m)$ contribute by  a rank two matrix.
\begin{lem}\label{w02} Let $n,m\in \Z$. Let $F(x)=q(U_n,U_m)(\log x)$, then one has :
\begin{equation}\label{hh}
W_{0,2}(V_n,V_m)=W_{0,2}^\#(F)=	\frac{32 L \sinh ^2\left(\frac{L}{4}\right) \left(L^2-16 \pi ^2 m n\right)}{\left(L^2+16 \pi ^2 m^2\right) \left(L^2+16 \pi ^2 n^2\right)}
\end{equation}
\end{lem}
 \proof This is best verified by direct computation. 
 \endproof 
 \subsection{The matrix $ W_p(V_n,V_m)$}
~The contribution of the non archimedean primes is given by \eqref{bombtestp}, \ie 
\begin{equation}\label{bomp}
\sum W_p(V_n,V_m)=\sum_{1<k\leq \exp(L)}\Lambda(k) k^{-1/2}q(U_n,U_m)(\log k).
 \end{equation}
 \subsection{The matrix $ W_\R(V_n,V_m)$}

 ~Let $\omega(x)=q(U_n,U_m)(x)$, then\eqref{bombtestR} gives
  \begin{align*}
W_\R(V_n,V_m) &=
\int_0^L\frac{\exp \left(\frac{x}{2}\right) \omega(x)- \omega(0)}{\exp \left(x\right)-\exp \left(-x\right)}dx- \omega(0)\int_L^{\infty}\frac{dx}{\exp \left(x\right)-\exp \left(-x\right)}\\ &+\frac{1}{2} (\gamma +\log (4 \pi )) \omega(0).
      \end{align*}
Since
$$
\int_L^{\infty}\frac{dx}{\exp \left(x\right)-\exp \left(-x\right)}=\frac 12 \log\left(\frac{e^L+1}{e^L-1} \right) 
$$
one obtains 
\begin{align}
W_\R(V_n,V_m) &=\frac{ \omega(0)}{2}\left(\gamma +\log \left(4 \pi \frac{e^L-1}{e^L+1}  \right)\right)\nonumber\\ +&\int_0^L\frac{\exp \left(\frac{x}{2}\right) \omega(x)- \omega(0)}{\exp \left(x\right)-\exp \left(-x\right)}dx.\label{weinfty}
\end{align}
The explicit expression for $\omega(x)=q(U_n,U_m)(x)$ for $x\in [0,L]$, is given by Lemma \ref{polarize0}.
The next step is to compute the integrals involved, they are given in terms of known functions in the following proposition \ref{computearch}. We let
$$
\rho(x):=\frac{\exp(x/2)}{\exp(x)-\exp(-x)}
$$
We let $\psi(z)=\Gamma^{\prime}(z) / \Gamma(z)$ be the digamma function and  $\psi ^{(1)}$ be its derivative. We use the standard notation for hypergeometric functions. We use the notation $\Phi$ for the Hurwitz-Lerch function 
$$\Phi (z,2,x)=\frac{1}{x^2}+\frac{z}{(x+1)^2}+\frac{z^2}{(x+2)^2}+\frac{z^3}{(x+3)^2}+\frac{z^4}{(x+4)^2}+\ldots $$
An  important feature of the formulas is that the parameter  $e^{-2L}$ in the various series involved is of modulus $<1$ thus ensuring the convergence and in fact fast numerical convergence for $L$ of order $10$.
\newpage
\begin{prop}\label{computearch} We use the symbols $\Re(z)$ and $\Im(z)$ for the real and imaginary parts of a complex number $z$. One has \begin{align}
	&\int_0^L\sin(2 \pi n x/L)\rho(x)dx=\\ & e^{-L/2}\Im(\frac{2L}{L+4 \pi in}\  _2F_1(1, \frac{\pi i n}{ L}+\frac{1}{4}; \frac{\pi i n}{ L}+\frac{5}{4};e^{-2L}))+\frac{1}{2}\Im(\psi(\frac{\pi i n}{ L}+\frac{1}{4})).\nonumber\\
	&\int_0^L x\ \cos(2 \pi n x/L)\rho(x)dx=\\ & -Le^{-L/2}\Im\left(\frac{2L}{4 \pi  n-i L}\, _2F_1\left(1,\frac{1}{4}+\frac{i n \pi }{L};\frac{5}{4}+\frac{i n \pi }{L};e^{-2L}\right)\right)\nonumber\\
	&-\frac{e^{-L/2}}{4}\Re\left(\Phi \left(e^{-2L},2,\frac{i \pi  n}{L}+\frac{1}{4}\right)\right)+\frac{1}{4}\Re\left(\psi ^{(1)}\left(\frac{\pi i n}{ L}+\frac{1}{4}\right)\right).\nonumber\\
	&\int_0^L \ (\cos(2 \pi n x/L)-1)\rho(x)dx=\\ &-e^{-L/2}\Re\left(\frac{2L}{L+4 \pi in}\, _2F_1\left(1,\frac{\pi i n}{ L}+\frac{1}{4};\frac{\pi i n}{ L}+\frac{5}{4};e^{-2L}\right)\right)\nonumber\\ \nonumber
	&+2e^{-L/2}\ _2F_1\left(\frac{1}{4},1;\frac{5}{4};e^{-2L}\right)-\frac{1}{2}\Re\left(\psi(\frac{\pi i n}{ L}+\frac{1}{4})-\psi(\frac 14)\right).
\end{align}	
\end{prop}
\proof In each case one first changes variables to $y=2\pi x/L$ and then expand, with $a:=\frac{L}{2\pi}$
$$
\rho(x)=\frac{\exp \left(\frac{a y}{2}\right) }{\exp (a y)-\exp (-a y)}=\sum_{k=0}^\infty \exp(b(k)y), \  \ b(k)=\frac{-a(1+4k) }{2}
$$
The change of variables introduces an overall factor $(\frac{L}{2\pi})^2$ for the middle integral and $(\frac{L}{2\pi})$ for the others. 
  One obtains in this way for each integral a sum of terms indexed by $k\in \N$ and which are
\begin{align}
\int_0^{2 \pi } \exp (b x) \sin (n x) \, dx=&\frac{n-e^{2 \pi  b} n}{b^2+n^2}\\
\int_0^{2 \pi } x \exp (b x) \cos (n x) \, dx=& \frac{2 \pi  e^{2 \pi  b}b-e^{2 \pi  b} +1}{b^2+n^2}-\frac{2 n^2\left(1-e^{2 \pi  b}\right)}{\left(b^2+n^2\right)^2}\\
\int_0^{2 \pi }  \exp (b x)( \cos (n x)-1) \, dx=&\frac{n^2-e^{2 \pi  b} n^2}{b^3+b n^2}	
\end{align}

\vspace{0.3cm}
All these expressions are affine in $e^{2 \pi  b}$ whose coefficient gives the general term of  a sum over $k\in \N$ using
$$
e^{2 \pi  b(k)}=\exp(-\frac L2(1+4k))=e^{-L/2} z^k, \ \ z=e^{-2L}.
$$ 
One then recognizes the series in $z$ involved and obtains the following expressions \begin{small}
\begin{align*}
&\int_0^L\sin(2 \pi n x/L)\rho(x)dx=\\
& e^{-L/2} \left(\frac{(i L) \, _2F_1\left(1,\frac{1}{4}-\frac{i n \pi }{L};\frac{5}{4}-\frac{i n \pi }{L};e^{-2L}\right)}{L-4 i \pi  n}-\frac{(i L) \, _2F_1\left(1,\frac{i \pi  n}{L}+\frac{1}{4};\frac{i \pi  n}{L}+\frac{5}{4};e^{-2L}\right)}{L+4 i \pi  n}\right)\\
&+\frac{1}{4} i \left(\psi ^{(0)}\left(\frac{1}{4}-\frac{i n \pi }{L}\right)-\psi ^{(0)}\left(\frac{i \pi  n}{L}+\frac{1}{4}\right)\right)
\end{align*}

\begin{align*}
&\int_0^L x\ \cos(2 \pi n x/L)\rho(x)dx=\\
&e^{-L/2} \left( 
    -\frac{1}{8} \, \Phi\left(e^{-2L}, 2, \frac{i \pi n}{L} + \frac{1}{4}\right)
    -\frac{1}{8} \, \Phi\left(e^{-2L}, 2, \frac{1}{4} - \frac{i \pi n}{L}\right)
\right) \\
& + e^{-L/2} \left(
    \frac{i L^2 \, {}_2F_1\left(1, \frac{i \pi n}{L} + \frac{1}{4}; \frac{i \pi n}{L} + \frac{5}{4}; e^{-2L}\right)}{4 \pi n - i L}
    - \frac{i L^2 \, {}_2F_1\left(1, \frac{1}{4} - \frac{i \pi n}{L}; \frac{5}{4} - \frac{i \pi n}{L}; e^{-2L}\right)}{4 \pi n + i L}
\right) \\
& + \frac{1}{8} \left(
    \psi^{(1)}\left(\frac{i \pi n}{L} + \frac{1}{4}\right) +
    \psi^{(1)}\left(\frac{1}{4} - \frac{i \pi n}{L}\right)
\right)\\
&\int_0^L \ (\cos(2 \pi n x/L)-1)\rho(x)dx=\\
&
 e^{-L/2} \left(-\frac{L \, _2F_1\left(1,\frac{i \pi  n}{L}+\frac{1}{4};\frac{i \pi  n}{L}+\frac{5}{4};e^{-2L}\right)}{L+4 i \pi  n}-\frac{L \, _2F_1\left(1,\frac{1}{4}-\frac{i n \pi }{L};\frac{5}{4}-\frac{i n \pi }{L};e^{-2L}\right)}{L-4 i \pi  n}\right)\\
&+2 e^{-L/2} \, _2F_1\left(\frac{1}{4},1;\frac{5}{4};e^{-2L}\right)+\frac{1}{4} \left(-\psi ^{(0)}\left(\frac{i \pi  n}{L}+\frac{1}{4}\right)-\psi ^{(0)}\left(\frac{1}{4}-\frac{i n \pi }{L}\right)+2 \psi ^{(0)}\left(\frac{1}{4}\right)\right)\end{align*}
\end{small}
Using the real and imaginary parts to simplify these expressions one obtains the required result. \endproof 
For the last integral we have computed a simplified form of the required expression which is 
\begin{equation}\label{corectc}
\int_0^L \ (\cos(2 \pi n x/L)-\exp(-x/2))\rho(x)dx=\int_0^L \ (\cos(2 \pi n x/L)-1)\rho(x)dx+c(L)\end{equation}
where the correction term is 
$$
c(L)=\int_0^L \frac{1-\exp \left(-\frac{x}{2}\right)}{\exp (x)-\exp (-x)} \, dx=$$ $$\log \left(e^{L/2}+1\right)+\frac{1}{4} \left(-2 \log \left(e^L+1\right)-\pi -\log (4)\right)+\tan ^{-1}\left(e^{L/2}\right)
$$
In fact we need to add the following to take into account the full Weil principal value
$$
w(L)=\frac{1}{2} (\gamma +\log (4 \pi ))-\frac{1}{2} \log \left(\frac{e^L+1}{e^L-1}\right)
$$
We then obtain 
$$
c(L)+w(L)=\frac{1}{2} \log \left(\frac{e^{L/2}-1}{e^{L/2}+1}\right)+\tan ^{-1}\left(e^{L/2}\right)-\frac{\pi }{4}+\frac{\gamma }{2}+\frac{1}{2} \log (8 \pi )
$$
To get lighter notations we let
\begin{align}\label{"funct}
&\alpha_L(n):=\frac{1}{\pi}\int_0^L\sin(2 \pi n x/L)\rho(x)dx, \ \\
&    \  \beta_L(n):=\frac 1L\int_0^L x\ \cos(2 \pi n x/L)\rho(x)dx, \  \\
& \gamma_L(n):=\int_0^L \ (\cos(2 \pi n x/L)-\exp(-x/2))\rho(x)dx+c(L)+w(L)
\end{align}
Using these notations, Proposition \ref{computearch} and \eqref{weinfty}, one obtains 
\begin{prop}\label{bigmatrix} The matrix $ W_\R(V_n,V_m)$ is given by the following table
 $$
 \begin{array}{ccc}
 	m\neq n &\vline & \frac{\alpha_L(m)-\alpha_L(n)}{n-m}\\
 	m= n &\vline & 2\gamma_L(n)-2\beta_L(n)
 \end{array}
$$
 \end{prop}
 \proof Follows from Lemma \ref{polarize0}.\endproof

\section{The infrared spectral triples}
In this section we shall construct infrared spectral triples naturally associated to the scaling operator in $L^2( [\lambda^{-1},\lambda],d^*u)$ with periodic boundary conditions. The intent is to modify the periodic boundary conditions in order to insert in the kernel of the modified scaling operator the eigenvector which realizes the minimum of the Weil quadratic form in the Hilbert space $L^2( [\lambda^{-1},\lambda],d^*u)$. In order to obtain this perturbation we work at the truncated level and use instead of the evaluation on the boundary of the interval $[\lambda^{-1},\lambda]$ an approximation to this evaluation which is given by the Dirichlet kernel. We then show the existence and uniqueness of the perturbed scaling operator, together with two fundamental facts. The first is that this operator becomes self-adjoint provided one modifies the inner product using the Weil quadratic form. The second point is that one can compute the spectrum of this operator using the Fourier transform of the minimal eigenvector. 
\subsection{Truncation of $QW_\lambda$}
 Let $\lambda>1$, $L=2\log \lambda$, $U_n$ as defined in \eqref{expbasis} and $V_n=\kappa(U_n)$ the orthonormal basis of $L^2( [\lambda^{-1},\lambda],d^*u)$ given in \eqref{vN}.
Let $N\in \N$, we consider the quadratic form $QW_\lambda^N$ obtained by restricting the Weil quadratic form $QW_\lambda$ to the finite dimensional space of test functions spanned by the functions $V_n$ for $\vert n\vert \leq N$.\newline
By \eqref{weilQexp}, the matrix elements of $T=QW_\lambda^N$ in the basis $V_n$ are given by
\begin{equation} \label{formN}
	\tau_{n,m}=\int_0^L q(U_n,U_m)(y)\cD(y)
  \qquad i,j\in \{-N,\ldots,N\}	
  \end{equation}
where $\cD$ is the real distribution $\cD=\log_*(\Psi^\sharp)$ on the interval $[0,L]$.
\begin{lem}\label{basicexpli} 
The matrix $\tau_{n,m}$ is a real symmetric matrix of the form
\begin{equation} \label{form}
	\tau_{i,i}=a_i,\quad  \forall i, \qquad\  \tau_{i,j}=\frac{b_i-b_j}{i-j}, \quad \forall j\neq i;
  \qquad i,j\in \{-N,\ldots,N\}	
  \end{equation}
where the real scalars $a_i$ fulfill $a_{-j}=a_j$ and $b_{-j}=-b_j$  for all $j\in \{-N,\ldots,N\}$.
\end{lem}
\proof This follows from the computation of the functions $q(U_n,U_m)(y)$ for $y\in [0,L]$ in \eqref{symrealmn} and \eqref{symrealn1}, which gives for $n\neq m$
$$
\tau_{n,m}=\int_0^L \frac{\sin(2\pi m y/L)-\sin(2\pi n y/L)}{\pi (n-m)} \cD(y)\implies $$ $$b_n=-\frac{1}{\pi}\int_0^L \sin(2\pi n y/L)
\cD(y)
$$
and for $n=m$,
$$
\tau_{n,n}=2\int_0^L(1-y/L) \cos(2\pi n y/L)\cD(y)=a_n
$$
\endproof 
\subsection{Properties of truncated matrices}
In this section we let $N$ be a positive integer, $E_N$ be the Hilbert space with orthonormal basis $\{V_n, n\in \{-N,\ldots,N\}\}$ and $T$ a real symmetric matrix of the form \eqref{form}. We recall the basic properties of \cite{CS} for matrices of this form. 
\begin{lem}\label{basics} $(i)$~Let $\gamma$ such that $\gamma(V_j):=V_{-j}$  $\forall j\in \{-N,\ldots,N\}$. One has $\gamma^2=\id$ and $T\gamma=\gamma T$.\newline
$(ii)$~Let $D$ be defined by $D(V_n):=n\,V_{n}$ for all $n\in \{-N,\ldots,N\}$. One has $D\gamma=-\gamma D$ and 
\begin{equation}\label{QD}
D\,	T-T\, D=\vert \beta\rangle\langle \eta\vert -\vert \eta\rangle\langle \beta\vert, \ \ \beta=\sum b_j\,V_j,\ \ \eta=\sum V_j.
\end{equation}	
\end{lem}
\proof $(i)$~One has $q_{-i,-j}=q_{i,j}$ for all $i,j\in \{-N,\ldots,N\}$.\newline
$(ii)$~The diagonal elements of the diagonal matrix $D$ are antisymmetric which gives $D\gamma=-\gamma D$. Let us prove \eqref{QD}. One has 
$$(DT)_{i,j}=i\tau_{i,j}, \ (TD)_{i,j}= j\tau_{i,j}$$ so that $(D\,	T-T\, D)_{i,j}=(b_i-b_j)$ for all $i,j\in \{-N,\ldots,N\}$. Similarly one has 
$$
(\vert \beta\rangle\langle \eta\vert)_{i,j}=\vert \beta\rangle_i\langle \eta\vert_j=b_i, \ \ (\vert \eta\rangle\langle \beta\vert)_{i,j}=\vert \eta\rangle_i\langle \beta\vert_j=b_j
$$
which gives the required equality.\endproof 
\begin{defn}\label{even-simple} A real symmetric matrix $T$ commuting with the $\Z/2$-grading $\gamma$ is {\bf even-simple} if its smallest eigenvalue is simple and the corresponding eigenvector $\xi$ satisfies $\gamma \xi=\xi$.
\end{defn}
We now assume that $T$ is even simple and positive and  let 
$\xi\in \ker\, T$, $\xi\neq 0$. The real symmetric positive  matrix $T$ defines an inner product on $\R^{2N+1}$ and its radical consists of the one dimensional subspace generated by $\xi$. Let us first show that we can normalize $\xi$ by the condition 
\begin{equation}\label{norxi}
	\langle \xi\mid \eta\rangle=1
\end{equation}
If $D\xi=0$ then $V_0\in \ker\, T$ fulfills \eqref{norxi}. So we can assume that $D\xi\neq 0$. One has 
 $TD\xi\neq 0$ since $D\xi$ is odd and linearly independent of $\xi$ while $\ker\, T$ is one-dimensional. By \eqref{QD} one has
 $$
 0\neq (D\,	T-T\, D)(\xi)=\vert \beta\rangle\langle \eta\vert\xi\rangle -\vert \eta\rangle\langle \beta\vert\xi\rangle=\vert \beta\rangle\langle \eta\vert\xi\rangle.
 $$ 
 Thus one can normalize $\xi$ so that $\langle \eta\vert\xi\rangle=1$. 
\begin{lem}\label{key} Assume $T\geq 0$ and $\Ker\, T=\C\xi$ where $\gamma \xi=\xi$ and $\langle \xi\mid \eta\rangle=1$. \newline
$(i)$~One has $T\, D\,\xi=-\beta $.\newline
$(ii)$~The operator $D':=D-\vert D\,\xi\rangle\langle \eta\vert$ induces a selfadjoint operator  $D"$ in the Hilbert space associated to the inner product defined by $T$, as quotient by null vectors.	\newline
$(iii)$~One has, denoting by $\xi_j$ the components of $\xi$, 
\begin{equation}
   \Det(D" - s) =   \Det(D - s) \sum_{j=-N}^N (j -s)^{-1} \xi_j.
   \label{eq:detDp}
   \end{equation}
\end{lem}
\proof $(i)$~We apply \eqref{QD} and get, using  $T\xi=0$ and $\langle \beta\vert \xi \rangle=0$ since the two eigenspaces of $\gamma$ are orthogonal,
$$
-T\, D\,\xi=(D\,	T-T\, D)\xi=\vert \beta\rangle\langle \eta\vert \xi\rangle -\vert \eta\rangle\langle \beta\vert \xi \rangle=\beta.
$$
$(ii)$~The inner product defined by $T$ is given by 
$$
\langle f\mid g\rangle_T=\langle T f\mid g\rangle.
$$
We first show that
\begin{equation}\label{herm}
	\langle D' f\mid g\rangle_T=\langle f\mid D'g\rangle_T, \  \forall f, g
\end{equation}
  One has, with $R=-\vert D\,\xi\rangle\langle \eta\vert$ 
$$
\langle D' f\mid g\rangle_T=\langle T D'f\mid g\rangle=\langle T Df\mid g\rangle+\langle T Rf\mid g\rangle.
$$
By $(i)$, one has $TR=-\vert TD\xi\rangle \langle \eta\vert=\vert \beta\rangle\langle \eta\vert$.  Thus
$$
T D'=TD+\vert \beta\rangle\langle \eta\vert
$$
Moreover by \eqref{QD}, one has $TD-DT=-\vert \beta\rangle\langle \eta\vert +\vert \eta\rangle\langle \beta\vert$. Thus
$$
T D'=DT+\vert \eta\rangle\langle \beta\vert,
$$ 
$$
\langle D' f\mid g\rangle_T=\langle DTf\mid g\rangle+\langle R'f\mid g\rangle, \ \ R'=\vert \eta\rangle\langle \beta\vert.
$$
 Moreover, using that both $T$ and $D$ are selfadjoint,  
 $$
 \langle f\mid D'g\rangle_T=\langle Tf\mid D g\rangle+\langle Tf\mid R g\rangle=\langle DTf\mid g\rangle+\langle f\mid TR g\rangle
 $$
 and the required equality follows from 
 $$
 TR g=(-TD\xi)\langle \eta\vert g\rangle=\beta \langle \eta\vert g\rangle
 $$
 $$
 \langle f\mid TR g\rangle=\langle f\mid (\vert \beta\rangle\langle \eta\vert g\rangle=\langle f\mid  \beta\rangle\langle \eta\vert g\rangle=\langle R'f\mid g\rangle.
 $$
 The Hilbert space $\cH$ obtained from $E_N$ using the inner product $\langle f\mid g\rangle_T$ is the quotient of $E_N$ by the radical $\Ker\, T=\C\xi$. By construction one has $D'\xi=0$ so that $D'$ induces an operator $D"$ in $\cH$ and $D"$ is selfadjoint by \eqref{herm}. \newline 
 $(iii)$~Let $v_j$ be an orthonormal basis of $\cH$ of eigenvectors for  $D"$ with eigenvalues $\lambda_j$. Let $w_j\in E_N$ be lifts of the $v_j$. One has $D"(v_j)=\lambda_jv_j$ and hence $D'(w_j)=\lambda_jw_j+s_j \xi$ for some  scalars $s_j$. Thus in the basis of $E_N$ formed by $\xi$ and the $w_j$, the matrix of $D'$ is triangular, with $0$ and the $\lambda_j$ on the diagonal. Thus  one gets 
 \begin{equation}\label{det0}
 {\rm Det}(D'-s)=-s\prod (\lambda_j-s)=-s\, {\rm Det}(D"-s)
\end{equation}
 We now compute ${\rm Det}(D'-s)$. We start by writing, in terms of $R =- \langle D \xi \rangle \langle \eta |$:
 $$
D' -s = D+ R-s = (D-s) \left (\id + (D-s)^{-1} R\right)
$$
Consequently 
$$
\Det(D' - s) = \Det(D-s) \Det (\id + (D-s)^{-1} R).
$$
To compute the second determinant we use the identity 
$$
 \Det ({\rm id} + A )=\sum_{k=0}^{\infty} {\rm Tr}\left(\wedge^k A\right)
$$
applied to the rank one operator $A= (D-s)^{-1} R$. The  higher exterior powers $\wedge^k A$ vanish for $k>1$ thus
$$
 \Det ({\rm id} + (D-s)^{-1} R )= 1 - {\rm Tr} \left( | (D-s)^{-1} D \xi \rangle \langle \eta | \right) = -s \langle \eta | (D-s)^{-1} \xi \rangle, 
 $$
 using $(D-s)^{-1} D \xi =\xi+ s (D-s)^{-1}\xi$ and $ \langle \eta | \xi\rangle=1$.
 Hence,
$$
   \Det(D' - s) = -s\, \Det (D-s) \langle \eta | (D-s)^{-1}\xi  \rangle = -s \, \Det (D-s)\, \sum_{j=-N}^N (j -s)^{-1} \xi_j.
   $$
   Thus one obtains \eqref{eq:detDp} using \eqref{det0}.
 \endproof 

\subsection{The Dirichlet Kernel $\delta_N$ as an approximation of the Dirac Delta}

The \textbf{Dirichlet kernel}  approximates the Dirac delta function as $N \to \infty$. It is \begin{equation}\label{dirichletdef}
D_N(x) = \sum_{n=-N}^{N} \exp(2 \pi i nx/L), \forall x\in [0, L]
\end{equation}
This can be simplified using the geometric series formula:
\begin{equation}
D_N(x) = \sin\left(\frac{\pi  (2 N+1) x}{L}\right)/\sin\left(\frac{\pi  x}{L}\right)
\end{equation}
for $x \neq 0$ (mod $L$) while $D_N(0)=D_N(L)=2N+1$.
\begin{lem}\label{dirichlet} Let $D=-i\partial_L$ be the selfadjoint operator of differentiation on  $L^2([0,L],dx)$ with periodic boundary conditions. Let $f\in \operatorname{Dom}D$. Then 
\begin{equation}
\lim_{N \to \infty} \frac{1}{L}\int_{0}^{L} D_N(x) f(x) \, dx = f(0)
\end{equation}	
\end{lem}
\proof Let $\hat{f}(n) = \frac{1}{L}\int_{0}^{L} f(x) e^{-2\pi inx/L} \, dx$ be the Fourier coefficients of $f$. 
We have by Parseval's theorem, with  a  constant $c>0$:
\begin{equation}
\sum_{n=-\infty}^{\infty} |n|^2 |\hat{f}(n)|^2 =c\, \|f'\|_{L^2}^2 < \infty
\end{equation}
Since $D_N(x) = \sum_{n=-N}^{N} e^{2\pi inx/L}$, we have:
\begin{equation}
\frac{1}{L}\int_{0}^{L} D_N(x) f(x) \, dx = \sum_{n=-N}^{N} \hat{f}(n)
\end{equation}
So we need to prove that the Fourier series of $f$ converges at $x=0$:
\begin{equation}
\sum_{n=-N}^{N} \hat{f}(n) \to f(0) \quad \text{as } N \to \infty
\end{equation}
 By the Cauchy-Schwarz inequality:
\begin{align}
\sum_{n \neq 0} |\hat{f}(n)| &= \sum_{n \neq 0} \frac{1}{|n|} \cdot |n||\hat{f}(n)|\nonumber \\
&\leq \left(\sum_{n \neq 0} \frac{1}{n^2}\right)^{1/2} \left(\sum_{n \neq 0} n^2|\hat{f}(n)|^2\right)^{1/2}< \infty\nonumber 
\end{align}
Therefore the Fourier series $\sum_{n=-\infty}^{\infty} \hat{f}(n)e^{2\pi inx/L}$ converges absolutely and uniformly  to $f(x)$ for all $x$, hence:
$$
\sum_{n=-N}^{N} \hat{f}(n) \to f(0) \quad \text{as } N \to \infty
$$
\endproof
The scaling operator is defined as
\begin{equation}
D_{\log}^{(\lambda)} =-i u\frac{\partial}{\partial u} =-i \frac{\partial}{\partial \log u}
\end{equation}
acting on  $L^2([\lambda^{-1},\lambda],d^*u)$ subject to periodic boundary conditions.
\begin{cor}\label{dirichlet1} Let $V_n(u):=U_n(\log (\lambda u))$ as in \eqref{vN}, $D_{\log}^{(\lambda)}$ be the scaling operator with periodic boundary conditions in  $L^2( [\lambda^{-1},\lambda],d^*u)$ and $f\in \operatorname{Dom}D_{\log}^{(\lambda)}$.
\begin{equation}
\lim_{N \to \infty} \langle \delta_N\mid f\rangle  = f(\lambda), \  \  \delta_N:=\frac{1}{\sqrt L}\  \sum_{n=-N}^{N}V_n 
\end{equation}	
\end{cor}
\proof~This follows from Lemma \ref{dirichlet} using the isometry $\kappa$ of \eqref{toa} to pass from $L^2([\lambda^{-1},\lambda],d^*u)$ to $L^2([0,L],dx)$, and the equality 
$$
D_N=\frac{1}{\sqrt L}\  \sum_{n=-N}^{N}U_n 
$$
which follows from \eqref{dirichletdef} and \eqref{expbasis}.\endproof 
\subsection{The perturbed scaling operator}

The perturbed scaling operator $\dln$ is obtained from the following :
\begin{prop}\label{pertscal} Let $\lambda>1$ and $N$ such that the truncated Weil quadratic form is even simple. Let $\xi$ be the corresponding even eigenvector. There exists a unique operator $\dln$ with the same domain as $D_{\log}^{(\lambda)}$ which agrees with this operator on the kernel of $\delta_N$ and such that $\dln(\xi)=0$.	
\end{prop}
\proof The Hilbert space  $L^2([\lambda^{-1},\lambda],d^*u)$ is the direct sum of the finite dimensional subspace $E_N$ spanned by the $V_n$ for $\vert n\vert \leq N$ and its orthogonal complement $E_N^\perp$. Since the $V_n$ form an orthonormal basis of eigenvectors for $D_{\log}^{(\lambda)}$, this operator splits as the direct sum of its restrictions to $E_N$ and $E_N^\perp$. The linear form $\delta_N$ vanishes on $E_N^\perp$ and $\xi \in E_N$ by construction. Thus the existence and uniqueness of $\dln$ is reduced to the finite dimensional subspace $E_N$. One has $\delta_N(\xi)\neq 0$ so $\Ker \delta_N$ and $\C\xi$ span $E_N$ and the equality 
$$
\dln(\beta)= D_{\log}^{(\lambda)}(\alpha), \ \forall \beta=\alpha +x\xi, \ \alpha \in \Ker \delta_N, \ x\in \C
$$ 
uniquely determines the operator $\dln$.\endproof

\subsection{Regularized Determinant}
The regularized determinant is defined by
\begin{equation}
\detreg(D - s) = \exp\left(-\zeta_D'(0; s)\right)
\end{equation}
where $\zeta_D(z; s):=\sum (\lambda-s)^{-z}$ is the spectral zeta function. There is an ambiguity in the definition of the zeta function (see \cite{SW}, \S 7.1) since raising $\lambda-s$ to the power $-z$ implies the choice of a determination of $\log(\lambda-s)$. There is a clear such choice for eigenvalues $\lambda\to +\infty$ but in the negative direction one needs to make a choice for $(-1)^{-z}$ and this choice of the spectral cut affects the result in the following basic example. Note that, in this example and taking for simplicity $L=2 \pi$, one could guess the regularized determinant from the Euler product as given by the sine function $\sin(\pi s)$ but this would violate the spectral invariance $s\to s+1$. The phase factor in front of the sine function repairs this violation.
\begin{lem} Let $L>0$. The regularized determinant for the Dirac operator $D$ with spectrum $\frac{2\pi}{L}\mathbb{Z}$ is given, using $(-1)^{-z}:=e^{-i\pi z}$, by 
	\begin{equation}\label{detreg0}
\detreg(D - s) = 1-e^{-iL s}
\end{equation}
\end{lem}
\proof 
We first take $L=2 \pi$ to simplify notations. The spectral zeta function is given by
\begin{equation}
\zeta_D(z; s) = \sum_{n \in \mathbb{Z} } (n-s)^{-z}
\end{equation}
This converges for $\mathrm{Re}(z) > 1$. One has 
\begin{equation}
\zeta_D(z; s) = \sum_{n=1}^{\infty} (n-s)^{-z} + \sum_{n=-\infty}^{0} (n-s)^{-z}
\end{equation}
For the second sum, substitute $n \to -m$ with $m \geq 0$ and use $(-1)^{-z}:=e^{-i\pi z}$ 
\begin{equation}
\sum_{n=-\infty}^{0} (n-s)^{-z} = \sum_{m=0}^{\infty} (-m-s)^{-z} = e^{-i\pi z} \sum_{m=0}^{\infty} (m+s)^{-z}
\end{equation}
 So that using the Hurwitz zeta function:
\begin{equation}
\zeta(z, a) = \sum_{n=0}^{\infty} (n+a)^{-z}
\end{equation}
one obtains 
$$
\zeta_D(z; s) = \sum_{n=0}^{\infty} (n+1-s)^{-z} + e^{-i\pi z} \sum_{m=0}^{\infty} (m+s)^{-z}= \zeta(z, 1-s)+ e^{-i\pi z}\zeta(z, s)$$

Thus one gets the equality 
\begin{equation}
\zeta_D'(0; s) =  \zeta'(0, 1-s)-i\pi \zeta(0, s)+ \zeta'(0, s)\end{equation}
Moreover one has the classical expressions
\begin{equation}
\zeta(0, a)=\frac{1}{2}-a,\ \zeta^{\prime}(0, a)=\log \Gamma(a)-\frac{1}{2} \log (2 \pi)
\end{equation}
which give 
\begin{equation}
\zeta_D'(0; s) =- \log (2 \pi) +\log \Gamma(s)+\log \Gamma(1-s) -i\pi(\frac{1}{2}-s)
\end{equation}
$$
\exp\left(-\zeta_D'(0; s)\right)=\frac{2\pi}{\Gamma(s)\Gamma(1-s)}\exp\left(i\pi(\frac{1}{2}-s)\right)=$$ $$=2 i \sin(\pi s)e^{-i\pi s}=1-e^{-2i\pi s}
$$
The general case follows using $\zeta_D(0; s)=0$. One has indeed, for $a>0$ 
$$
\zeta_{a D}(z;a s)=a^{-z}\zeta_D(z; s)\implies \zeta_{a D}'(0;a s)=-\log a\ \zeta_D(0; s)+\zeta_D'(0; s)
$$
which gives the required result for $a=2\pi/L$.
\endproof 

 \subsection{Spectrum and regularized determinant of $\dln$}
 In this section we prove the fundamental properties of $\dln$ and show that its spectrum is real and that its regularized determinant is (up to a phase factor) the Fourier transform of the minimal eigenvector $\xi$.\newline
We first compute the Fourier transform of functions on $[\lambda^{-1},\lambda]$  extended by $0$ to  $\R_+^*$. The Fourier transform for the duality $\langle \R_+^*\mid \R\rangle$ is defined by
$$
\fourier_\mu (f)(s):=\int_{\R_+^*} f(u)u^{-is}d^*u.
$$

\begin{prop}  Let $\lambda>1$, $L=2\log \lambda$, $N\in \N$, $\xi_j\in \C$ for $j\in \{-N,\ldots,N\}$, $$
\xi(u):=\sum_{\{-N,\ldots,N\}} \xi_k\, V_k(u), \ \forall u\in [\lambda^{-1},\lambda], \ \xi(u)=0, \ \forall u\notin [\lambda^{-1},\lambda]. 
$$ The Fourier transform $\widehat\xi=\fourier_\mu(\xi)$ is the entire function given by
\begin{equation}\label{four}
\widehat\xi(z)=2\,L^{-1/2}\sin \left(zL/2\right)\left(\sum_{\{-N,\ldots,N\}} \frac{\xi_j}{z-2\pi j/L}\right).
\end{equation}	
\end{prop}
\proof One has for $k\in \Z$, using $x=\log(\lambda u)$, $u^{-is}=\lambda^{is}\exp (-i s x)$, $d^*u=dx$,
$$
\int_{\lambda^{-1}}^\lambda V_k(u)u^{-is}d^*u=\int_{\lambda^{-1}}^\lambda U_k(\log(\lambda u))u^{-is}d^*u=
$$
	$$
	=L^{-1/2}\lambda^{is}\int_0^L \exp (2 \pi  i k x/L) \exp (-i s x) \, dx=-i\frac{1 -e^{-isL} }{s-2 \pi  k/L}L^{-1/2}e^{isL/2}.
	$$
	Thus the  Fourier transform of $\xi(u)$ is 
	$$
	\widehat\xi(z)=2\,L^{-1/2}\sin \left(zL/2\right)\left(\sum_{\{-N,\ldots,N\}} \frac{\xi_j}{z-2\pi j/L}\right).
	$$
	The zeros $z\in 2\pi \Z/L$ of $\sin \left(sL/2\right)$ cancell the poles at $2\pi j/L$ which occur when $\xi_j\neq 0$ and remain as zeros of $\widehat\xi(z)$ otherwise. \endproof 
 \begin{theorem}\label{finmain} Let $\epsilon_N$ be the smallest eigenvalue of $QW_\lambda^N$  assumed simple and  $\xi$ the corresponding eigenvector assumed even, normalized by $\delta_N(\xi)=1$.\newline
 $(i)$~The operator $\dln$ is selfadjoint in the direct sum $E'_N\oplus E_N^\perp$ where  on the subspace $E'_N=E_N/\C\xi$ the inner product is given 	by the restriction of the quadratic form $QW_\lambda^N-\epsilon_N \langle \mid \rangle$.\newline 
 $(ii)$~The regularized determinant of $\dln$ is given by $$\detreg(\dln - z)=-i\,\lambda^{-iz}\widehat\xi(z)$$ where $\widehat\xi$ is the Fourier transform of $\xi$ for the duality $\langle\R_+^*\mid \R\rangle$.\newline
$(iii)$~The Fourier transform $\widehat\xi(z)$ is an entire function, all its zeros are on the real line and coincide with the spectrum of $\dln$.
 \end{theorem}
 \proof $(i)$~We apply Lemma \ref{key} to $T:=QW_\lambda^N-\epsilon_N\id
 $ where $\epsilon_N$ is the smallest eigenvalue of $QW_\lambda^N$ which is assumed to be simple and even. Let $D$, $D'$ and $D"$ be the associated operators. One has 
 \begin{equation}\label{propor}
 	D_{\log}^{(\lambda)}\vert_{E_N}=\frac{2\pi}{L}D, \ \ \dln\vert_{E_N}=\frac{2\pi}{L}D', \ \ \dln\vert_{E'_N}=\frac{2\pi}{L}D"
 \end{equation}
 By construction the operator $\dln$ decomposes the direct sum in the direct sum $E'_N\oplus E_N^\perp$. Hence the result follows from Lemma \ref{key}, $(ii)$.
 \newline 
 $(ii)$~By \eqref{propor}  one has 
 $$
 \Det(D_{\log}^{(\lambda)}\vert_{E_N}-\frac{2\pi}{L} s)=(\frac{2\pi}{L})^{2N+1}\Det(D - s)
 $$
 $$
 \Det(\dln\vert_{E'_N} -\frac{2\pi}{L} s)=(\frac{2\pi}{L})^{2N}\Det(D" - s) 
 $$ 
 We now apply Lemma \ref{key}, $(iii)$, but note that the normalization of $\xi$ given by $\delta_N(\xi)=1$ differs from $\langle \xi'\mid \eta\rangle=1$. One has $\delta_N=L^{-1/2}\eta$ by Corollary \ref{dirichlet1}  thus giving $\xi' =L^{-1/2}\xi $. Lemma \ref{key}, $(iii)$ shows that 
 $$
 \Det(\dln\vert_{E'_N} -\frac{2\pi}{L} s) =\frac{L}{2\pi } \Det(D_{\log}^{(\lambda)}\vert_{E_N}-\frac{2\pi}{L} s)\sum_{j=-N}^N (j -s)^{-1} \xi'_j
 $$
 which with $\xi' =L^{-1/2}\xi $,  $z=\frac{2\pi}{L} s$ gives 
 $$
 \Det(\dln\vert_{E'_N} -z) = \Det(D_{\log}^{(\lambda)}\vert_{E_N}-z)\frac{L^{1/2}}{2\pi }\sum_{j=-N}^N (j -\frac{L}{2\pi }z)^{-1} \xi_j
 $$
 It then follows from the multiplicativity of the regularized determinant that 
 \begin{equation}\label{propor1}
 \detreg(\dln - z)=L^{-1/2}\detreg(D_{\log}^{(\lambda)}-z)\left(\sum_{\{-N,\ldots,N\}} \frac{\xi_j}{2\pi j/L-z}\right)
  \end{equation}
By \eqref{four}, the  Fourier transform of $\xi(u)$ is  
	$$
	\widehat\xi(z)=2\,L^{-1/2}\sin \left(zL/2\right)\left(\sum_{\{-N,\ldots,N\}} \frac{\xi_j}{z-2\pi j/L}\right)
	$$
	By \eqref{detreg0} one has
	$$
	\detreg(D_{\log}^{(\lambda)}-z)=1-\exp(-iz\, L)=2i\,\exp(-iz\, L/2)\sin \left(zL/2\right),
	$$ 
	thus \eqref{propor1} gives 
	$$
	\detreg(\dln - z)=L^{-1/2}\left(1-\exp(-iz\, L)\right)\left(\sum_{\{-N,\ldots,N\}} \frac{\xi_j}{2\pi j/L-z}\right)=$$
	$$=-i\,\exp(-iz\, L/2)\,\widehat\xi(z)=-i\,\lambda^{-iz}\,\widehat\xi(z)
	$$
	$(iii)$~The Fourier transform $\widehat\xi$ is an entire function since $\xi$ is an $L^1$-function with compact support. The regularized determinant $\detreg(\dln - z)$ is the product of the determinants 
	$$
	\detreg(\dln - z)=\Det(\dln\vert_{E'_N} -z) \detreg(D_{\log}^{(\lambda)}\vert_{E_N^\perp} - z)
	$$
	In this factorization the first term is the characteristic polynomial of a selfadjoint matrix and hence all its zeros are real. The zeros of the second term form the set 
	$$
	\{2\pi j/L\mid j\in \Z, \vert j\vert >N\}
	$$ 
	which gives the required result.
	  \endproof 
\section{Numerical results}\label{sectnumerical}
	One computes using the above formulas the matrix of the Weil quadratic form and the spectrum of the operator $\dln$. These computations require high precision but are easily preformed using 200 digits accuracy due to the fast convergence of the special functions involved. 	
 The first case is $\lambda=3$ and one takes $N=120$.
	\begin{figure}[H]
  \centering
\includegraphics[scale=0.8]{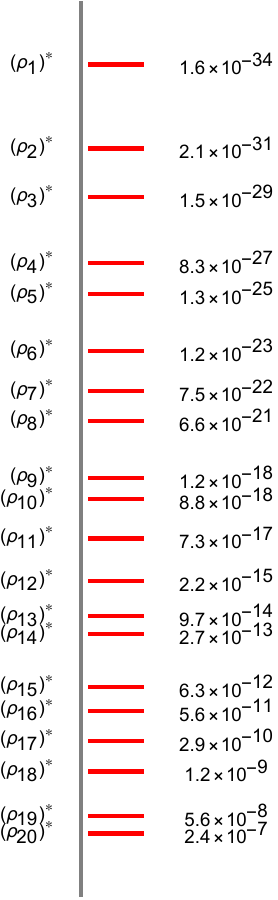}\label{toep}
\caption{This shows the differences between  the  first twenty zeros of $\zeta\left(\frac{1}{2}+i s\right)$  and the eigenvalues of the operator $\dln$ for $\lambda=3$ and $N=120$.}
\end{figure}
We then consider the first fifty zeros of zeta, use still $N=120$ and the values of $\lambda$ given by 
$$
\lambda=\sqrt 12\sim 3.4641, \ \lambda=\sqrt 13\sim 3.60555, \  \lambda =\sqrt 14 \sim 3.74166. \ 
$$
We get the following table giving an upper bound on the absolute value of the difference between the nontrivial zeros of the Riemann zeta function $\zeta\left(\frac{1}{2}+i s\right)$  and the eigenvalues of the operator $\dln$. 
\begin{center}	
\begin{tabular}{|l|l|l|l|}
\hline & $\lambda=\sqrt{12}$ & $\lambda=\sqrt{13}$ & $\lambda=\sqrt{14}$ \\
\hline 1 & $3.41 \times 10^{-50}$ & $2.44 \times 10^{-55}$ & $1.07 \times 10^{-60}$ \\
\hline 2 & $5.89 \times 10^{-47}$ & $4.5 \times 10^{-52}$ & $2.08 \times 10^{-57}$ \\
\hline 3 & $5.18 \times 10^{-45}$ & $4.16 \times 10^{-50}$ & $2 . \times 10^{-55}$ \\
\hline 4 & $4.2 \times 10^{-42}$ & $3.65 \times 10^{-47}$ & $1.89 \times 10^{-52}$ \\
\hline 5 & $7.84 \times 10^{-41}$ & $7.11 \times 10^{-46}$ & $3.81 \times 10^{-51}$ \\
\hline 6 & $1.07 \times 10^{-38}$ & $1.06 \times 10^{-43}$ & $6.13 \times 10^{-49}$ \\
\hline 7 & $9.42 \times 10^{-37}$ & $1 . \times 10^{-41}$ & $6.17 \times 10^{-47}$ \\
\hline 8 & $1.05 \times 10^{-35}$ & $1.19 \times 10^{-40}$ & $7.66 \times 10^{-46}$ \\
\hline 9 & $3.25 \times 10^{-33}$ & $4.12 \times 10^{-38}$ & $2.94 \times 10^{-43}$ \\
\hline 10 & $2.99 \times 10^{-32}$ & $3.98 \times 10^{-37}$ & $2.96 \times 10^{-42}$ \\
\hline 11 & $3.76 \times 10^{-31}$ & $5.5 \times 10^{-36}$ & $4.42 \times 10^{-41}$ \\
\hline 12 & $1.87 \times 10^{-29}$ & $3.04 \times 10^{-34}$ & $2.68 \times 10^{-39}$ \\
\hline 13 & $1.28 \times 10^{-27}$ & $2.29 \times 10^{-32}$ & $2.19 \times 10^{-37}$ \\
\hline 14 & $4.47 \times 10^{-27}$ & $8.46 \times 10^{-32}$ & $8.43 \times 10^{-37}$ \\
\hline 15 & $2.18 \times 10^{-25}$ & $4.82 \times 10^{-30}$ & $5.48 \times 10^{-35}$ \\
\hline 16 & $2.76 \times 10^{-24}$ & $6.61 \times 10^{-29}$ & $8.02 \times 10^{-34}$ \\
\hline 17 & $2.3 \times 10^{-23}$ & $6.08 \times 10^{-28}$ & $8.02 \times 10^{-33}$ \\
\hline 18 & $1.59 \times 10^{-22}$ & $4.66 \times 10^{-27}$ & $6.73 \times 10^{-32}$ \\
\hline 19 & $1.59 \times 10^{-20}$ & $5.5 \times 10^{-25}$ & $9.11 \times 10^{-30}$ \\
\hline 20 & $9.55 \times 10^{-20}$ & $3.54 \times 10^{-24}$ & $6.19 \times 10^{-29}$ \\
\hline 21 & $2.36 \times 10^{-19}$ & $9.75 \times 10^{-24}$ & $1.86 \times 10^{-28}$ \\
\hline 22 & $6.7 \times 10^{-18}$ & $3.31 \times 10^{-22}$ & $7.38 \times 10^{-27}$ \\
\hline 23 & $5.24 \times 10^{-17}$ & $2.86 \times 10^{-21}$ & $6.89 \times 10^{-26}$ \\
\hline 24 & $8.4 \times 10^{-16}$ & $5.32 \times 10^{-20}$ & $1.45 \times 10^{-24}$ \\
\hline
 25 & $1.94 \times 10^{-15}$ & $1.33 \times 10^{-19}$ & $3.89 \times 10^{-24}$ \\
\hline 26 & $2.42 \times 10^{-14}$ & $2.07 \times 10^{-18}$ & $7.23 \times 10^{-23}$ \\
\hline 27 & $6.05 \times 10^{-13}$ & $5.94 \times 10^{-17}$ & $2.33 \times 10^{-21}$ \\
\hline 28 & $1.26 \times 10^{-12}$ & $1.34 \times 10^{-16}$ & $5.58 \times 10^{-21}$ \\
\hline 29 & $3.15 \times 10^{-12}$ & $4.09 \times 10^{-16}$ & $2.01 \times 10^{-20}$ \\
\hline 30 & $2.72 \times 10^{-11}$ & $4.21 \times 10^{-15}$ & $2.39 \times 10^{-19}$ \\
\hline 31 & $3.57 \times 10^{-10}$ & $6.61 \times 10^{-14}$ & $4.33 \times 10^{-18}$ \\
\hline 32 & $1.7 \times 10^{-9}$ & $3.6 \times 10^{-13}$ & $2.62 \times 10^{-17}$ \\
\hline 33 & $2.33 \times 10^{-9}$ & $5.66 \times 10^{-13}$ & $4.6 \times 10^{-17}$ \\
\hline 34 & $1.2 \times 10^{-7}$ & $4.03 \times 10^{-11}$ & $4.23 \times 10^{-15}$ \\
\hline 35 & $2.89 \times 10^{-7}$ & $1.04 \times 10^{-10}$ & $1.16 \times 10^{-14}$ \\
\hline 36 & $4.1 \times 10^{-7}$ & $1.85 \times 10^{-10}$ & $2.45 \times 10^{-14}$ \\
\hline 37 & $9.11 \times 10^{-7}$ & $4.92 \times 10^{-10}$ & $7.53 \times 10^{-14}$ \\
\hline 38 & $2.78 \times 10^{-6}$ & $1.94 \times 10^{-9}$ & $3.61 \times 10^{-13}$ \\
\hline 39 & $3.53 \times 10^{-5}$ & $3.24 \times 10^{-8}$ & $7.44 \times 10^{-12}$ \\
\hline 40 & $1.83 \times 10^{-4}$ & $2 . \times 10^{-7}$ & $5.24 \times 10^{-11}$ \\
\hline
\end{tabular}	
\end{center}
\newpage
\begin{center}
\begin{tabular}{|l|l|l|l|}
\hline 40 & $1.83 \times 10^{-4}$ & 2. $\times 10^{-7}$ & $5.24 \times 10^{-11}$ \\
\hline 41 & $1.67 \times 10^{-4}$ & $2.12 \times 10^{-7}$ & $6.22 \times 10^{-11}$ \\
\hline 42 & $2.97 \times 10^{-4}$ & $5.66 \times 10^{-7}$ & $2.23 \times 10^{-10}$ \\
\hline 43 & $2.19 \times 10^{-3}$ & $5.49 \times 10^{-6}$ & $2.64 \times 10^{-9}$ \\
\hline 44 & $4.35 \times 10^{-3}$ & $1.35 \times 10^{-5}$ & $7.51 \times 10^{-9}$ \\
\hline 45 & $1.19 \times 10^{-2}$ & $5.3 \times 10^{-5}$ & $3.8 \times 10^{-8}$ \\
\hline 46 & $1.27 \times 10^{-2}$ & $6.88 \times 10^{-5}$ & $5.65 \times 10^{-8}$ \\
\hline 47 & $2.87 \times 10^{-2}$ & $3.01 \times 10^{-4}$ & $3.66 \times 10^{-7}$ \\
\hline 48 & $1.43 \times 10^{-1}$ & $2 . \times 10^{-3}$ & $2.98 \times 10^{-6}$ \\
\hline 49 & $1.98 \times 10^{-1}$ & $3.01 \times 10^{-3}$ & $5.34 \times 10^{-6}$ \\
\hline 50 & $9.02 \times 10^{-2}$ & $2.04 \times 10^{-3}$ & $4.78 \times 10^{-6}$ \\
\hline
\end{tabular}
\end{center}
\section{Outlook}\label{outlook}

These numerical results provide evidence that the spectra of the operators $\dln$ tend\footnote{when $N\to\infty$ and $\lambda \to \infty$} to the  nontrivial zeros of the Riemann zeta function $\zeta(\frac 12+is)$. Establishing this convergence rigorously would amount to a proof of the Riemann Hypothesis. \newline
One can be even more ambitious using Theorem \ref{finmain} together with the observation of \cite{VJ} that the eigenfunction associated with the lowest eigenvalue of $QW_\lambda$ is well approximated by prolate spheroidal wave functions. As we explain now this suggests that  the regularized determinants 
$\detreg(\dln-s)$ behave as follows 
\begin{itemize}
\item For fixed $\lambda$, the functions $\detreg(\dln-s)$ converge\footnote{uniformly on compact subsets of $\C$} when $N\to \infty$ to the function $-i\,\lambda^{-iz}\widehat\xi_\lambda(z)$ where $\xi_\lambda$ is the eigenfunction of $QW_\lambda$ for the smallest eigenvalue, normalized by $\xi(\lambda)=1$.
\item When $\lambda\to \infty $ the functions $\widehat\xi_\lambda(z)$ multiplied by suitable constants, converge  uniformly on closed substrips of the open  strip $\Im(z)<\frac 12$ towards the $\Xi$-function of Riemann
$$
\Xi(s) = \xi(1/2 + is), \  \ \xi(z) = \frac{1}{2}z(z-1)\pi^{-z/2}\Gamma(z/2)\zeta(z)
$$  	
\end{itemize}
In other words the regularized determinants $\detreg(\dln-s)$ suitably multiplied by a factor of the form $e^{a+ib s}$ converge towards   $\Xi(s)$. This convergence would entail RH using Hurwitz theorem on the zeros of limits of holomorphic functions. \newline
We now give more details to give substance to this strategy. First when reading Riemann's paper \cite{Riemann} one finds that, using modern terminology, he understood his $\Xi$-function as the Fourier transform for the duality $\langle\R_+^*\mid \R\rangle$ of the function 
\begin{equation}\label{ku}	
k(u)=\cE(h)(u), \ \ h(u)=\frac \pi 2 u^2\left(2 \pi   u^2-3\right)e^{-\pi u^2}.
\end{equation}
where one uses the following map $\cE$:
\begin{equation}\label{emap}	
\cE(f)(u):=u^{1/2}\sum_1^\infty f(nu)
\end{equation}
The function $h(u)$ can be characterized as follows. One considers the Hermite operator (harmonic oscillator)
 \begin{equation}\label{herm}	\mathbf{H}f(u):=-f''(u)+4 \pi ^2 u^2 f(u)
 \end{equation} 
 and lets  $h_n$ be the normalized eigenfunction for the eigenvalue $2\pi(1+2n)$. These functions are even for $n$ even and, for $n$ multiple of $4$, invariant under the Fourier transform for the duality $\langle\R\mid \R\rangle$ which is defined by
\[
\mathbb{F}_{e_{\mathbb{R}}}(f)(y) := \int_{\mathbb{R}} f(x) e^{2\pi i x y} \, dx.
\]

\begin{lem}\label{hermfact}
The $\Xi$ function of Riemann is the Fourier transform of $k=\cE(h)$ where 	$h$ is, up to a multiplicative scalar, the only linear combination of $h_0, h_4$ with vanishing integral. More precisely one has, in terms of the normalized $h_n$
\begin{equation}\label{xih}	
  h=\frac{\sqrt{3}}{2^{11/4}}h_4-\frac{3}{2^{17/4}}h_0,  \ \  \text{and} \  \  \Vert h\Vert=\frac{\sqrt{33}}{2^{17/4}}.
\end{equation}
\end{lem}
\proof The normalized forms of $h_0, h_4$ are 
$$
h_0(x)=2^{1/4}e^{-\pi x^2}, \  \  h_4(x)=\left(\frac{16 \pi ^2 x^4-24 \pi  x^2+3}{2 \sqrt[4]{2} \sqrt{3}}\right) e^{-\pi x^2}
$$
Thus
$$
\frac{3}{2^{17/4}}h_0(x)=\frac{3}{16}e^{-\pi x^2}, \  \ \frac{\sqrt{3}}{2^{11/4}}h_4(x)=\left(\pi ^2 x^4-\frac{3 \pi  x^2}{2}+\frac{3}{16}\right)e^{-\pi x^2}
$$
which gives the required result using \eqref{ku}.\endproof
We  recall the construction of \cite{VJ} of an educated guess $k_\lambda$ for an approximation of a scalar multiple of  $\xi_\lambda$. It is based on the deformation of the harmonic oscillator called the prolate wave operator 
\begin{equation}\label{wop}
PW_\lambda := -\partial_x\left( (\lambda^2 - x^2) \partial_x \right) + (2\pi \lambda x)^2.
\end{equation} 
The eigenfunctions $h_{n,\lambda}(u)$ of $PW_\lambda$ have the same labelling as the Hermite functions $h_n$, they are even for $n$ even and invariant under the Fourier transform for $n$ multiple of $4$. 
In  agreement with Lemma \ref{hermfact}, the  educated guess $k_\lambda$ is \begin{equation}\label{ktoh}
  k_\lambda(u):=\cE(h_\lambda)(u), \ \ \forall u\in [\lambda^{-1},\lambda]
  \end{equation}
   where $h_\lambda$ is, up to a multiplicative scalar, the only linear combination of $h_{0,\lambda}, h_{4,\lambda}$ with vanishing integral. We refer to \cite{VJ}, Section 3, for the motivation behind the formula for $k_\lambda$ and the numerical evidence showing that it gives an approximation of a scalar multiple of  $\xi_\lambda$. Justifying rigorously this step is the main remaining obstacle to our approach to RH.\newline
We now use the educated guess \eqref{ktoh} and evaluate its convergence in the next lemma \ref{hermfact1}. We shall first describe  an estimate from \cite{MS}.
 \begin{lem}\label{meixnerlem}$(i)$~The eigenfunctions $h_{n,\lambda}$ of $PW_\lambda$, suitably normalized, fulfill for $n=0,4$  an estimate of the form (with $c<\infty$)
 \begin{equation}\label{esti01}
 \max_{x\in [-\lambda,\lambda]}\vert h_{n,\lambda}(x)-h_n(x)\vert\leq c\, \lambda^{-2}
\end{equation}
 $(ii)$~Let $h_\lambda$ be  the suitably normalized linear combination of $h_{0,\lambda}, h_{4,\lambda}$ with vanishing integral.
 	One has an estimate of the form (with $c<\infty$)
\begin{equation}\label{esti1}\max_{x\in [-\lambda,\lambda]}\vert h_\lambda(x)-h(x)\vert\leq c\, \lambda^{-2}
\end{equation}
 \end{lem}
\proof $(i)$~This follows from \cite{MS}, Satz 9, page 243, Section 3.2. entitled  "Die Sphäroidfunktionen $\operatorname{ps}_n^m\left(z ; \gamma^2\right)$" which asserts that uniformly for $z\in [-1,1]$ one has the estimate 
$$
\operatorname{ps}_n^m\left(z ; \gamma^2\right)=(-1)^m\left(\frac{4 \gamma}{\pi}\right)^{\frac{1}{4}} \frac{1}{(n-m)!}\left(\frac{(n+m)!}{2 n+1}\right)^{\frac{1}{2}}\left(1-z^2\right)^{m/ 2} D_{n-m}\left((2 \gamma)^{\frac{1}{2}} z\right)+O\left(\gamma^{-\frac{3}{4}}\right)$$
We need to explain carefully the notations of \cite{MS}. The differential equation defining the prolate spheroidal functions uses the operator  $$F_\gamma y \equiv \frac{d}{d z}\left[\left(1-z^2\right) \frac{d y}{d z}\right]+\left[\frac{-m^2}{1-z^2}+\gamma^2\left(1-z^2\right)\right] y$$
and we are only interested in the case when the angular parameter $m=0$. In that case the operator simplifies to
$$F_\gamma y \equiv \frac{d}{d z}\left[\left(1-z^2\right) \frac{d y}{d z}\right]+\gamma^2\,\left(1-z^2\right) y$$ 
We first relate this operator (for $m=0$) to the prolate wave operator of \eqref{wop}.
One lets $z=x/\lambda$ which gives 
$$
\left(\frac{d}{d z}\right)^2\mapsto \lambda^2\left(\partial_x\right)^2, \ \ -\gamma^2\,z^2\mapsto -\gamma^2/\lambda^2\,x^2
$$
Thus we need an overall minus sign and also 
\begin{equation}\label{wop1}
\gamma^2/\lambda^2=4\pi^2 \lambda^2\implies \gamma=2 \pi \lambda^2
\end{equation} 
The prolate spherical functions  $\operatorname{ps}_n:=\operatorname{ps}_n^0$  are related to those in Mathematica by 
$$
\operatorname{ps}_n\left(s ; \gamma^2\right)=PS(n,0,\gamma,s)
$$
The statement in \cite{MS} Satz 9 simplifies for $m=0$ to the estimate, uniform for $z\in [-1,1]$,
$$
\operatorname{ps}_n\left(z ; \gamma^2\right)=\left(\frac{4 \gamma}{\pi}\right)^{\frac{1}{4}}  \left(\frac{1}{(2 n+1)n!}\right)^{\frac{1}{2}} D_{n}\left((2 \gamma)^{\frac{1}{2}} z\right)+O\left(\gamma^{-\frac{3}{4}}\right)$$
which gives uniformly in the interval $[-\lambda,\lambda]$ 
\begin{equation}\label{meix0}
\left(\frac{4 \gamma}{\pi}\right)^{-\frac{1}{4}}\operatorname{ps}_n\left(x/\lambda ; \gamma^2\right)=  \left(\frac{1}{(2 n+1)n!}\right)^{\frac{1}{2}} D_{n}\left((2 \gamma)^{\frac{1}{2}} x/\lambda\right)+O\left(\gamma^{-1}\right)
 \end{equation}
 The hermite functions $D_n$ are defined in \cite{MS} by the equation 
 $$D_p(x)=(-1)^p e^{\frac{x^2}{4}} \frac{d^p}{d x^p} e^{-\frac{x^2}{2}}, \  \ D_0(x)=e^{-\frac{x^2}{4}} $$
  In our case we have  $\gamma=2\pi \lambda^2$, which gives
  $$
  (2 \gamma)^{\frac{1}{2}}/\lambda=(4\pi)^{\frac{1}{2}}\implies  D_{n}\left((2 \gamma)^{\frac{1}{2}} x/\lambda\right)=D_n\left((4\pi)^{\frac{1}{2}}x\right)
  $$
 One has in particular 
$$
 	D_0\left((4\pi)^{\frac{1}{2}}x\right)=e^{-\pi x^2}, \ \ D_4\left((4\pi)^{\frac{1}{2}}x\right)=e^{-\pi  x^2} \left(16 \pi ^2 x^4-24 \pi  x^2+3\right)
$$
 Thus in terms of the normalized Hermite functions $h_n$ we have 
 \begin{equation}\label{meix1}
D_0\left((4\pi)^{\frac{1}{2}}x\right)=2^{-1/4}h_0(x), \  \ 	D_4\left((4\pi)^{\frac{1}{2}}x\right)=2^{5/4}  \sqrt{3}\,h_4(x)
 \end{equation}
 We can thus use \eqref{meix0} to normalize the products $h_{n,\lambda}=c_n\lambda^{-1/2}\operatorname{ps}_n\left(x/\lambda ; \gamma^2\right)$  so that 
\begin{equation}\label{esti}
\max_{x\in [-\lambda,\lambda]}\vert h_{n,\lambda}(x)-h_n(x)\vert\leq c\, \lambda^{-2}, \  \ n=0,4.
\end{equation}
$(ii)$~By $(i)$, we control the values $h_{n,\lambda}(0)$ for $n=0,4$. The fundamental property of the prolate wave functions  is that they are eigenfunctions   of the compression of the Fourier transform by the orthogonal projection $P_\lambda$ of $L^2(\R)^{\rm even}$ on the subspace of functions with support in the  interval  $[-\lambda,\lambda]$.
Moreover for small values of $n$ such as $n=0,4$ the eigenvalues $\chi(\lambda)$ are such that  $1-\chi(\lambda)$ decays extremely fast when $\lambda\to \infty$.  For instance for $n=4$, by \cite{Fuchs}, Theorem 1, one has
$$1-\chi(\lambda)\sim\frac{2^{14}}{3} \sqrt{2} \pi ^5 e^{-4 \pi  \lambda^2+9 \log(\lambda)} $$
Moreover we have 
$$
\int h_{n,\lambda}(x)dx=\widehat{h_{n,\lambda}}(0)=\chi_n(\lambda)h_{n,\lambda}(0)
$$
By \eqref{esti} we control the differences $\vert h_{n,\lambda}(0)-h_n(0)\vert$ and hence the differences 
$$
\vert \int h_{n,\lambda}(x)dx-h_n(0)\vert=O(\lambda^{-2})
$$
It follows that we obtain a linear combination $h_\lambda$ of $h_{0,\lambda}$ and  $h_{4,\lambda}$ which has vanishing integral and fulfills \eqref{esti1}. \endproof 
In fact we show in the following two Figures the behavior of the functions 
$$
e_n(\lambda^2):=\lambda^2 \max_{x\in [-\lambda,\lambda]}\vert h_{n,\lambda}(x)-h_n(x)\vert
 $$
 \begin{figure}[H]
\begin{center}
\includegraphics[scale=0.64]{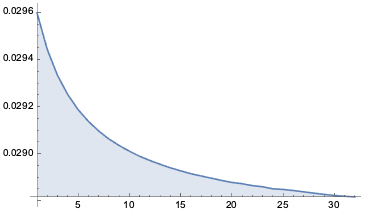}\\
\caption{Graph of $e_0(\mu)$ for $\mu\leq 36$.\label{fpro3}}
\end{center}
\end{figure}
\begin{figure}[H]
\begin{center}
\includegraphics[scale=0.64]{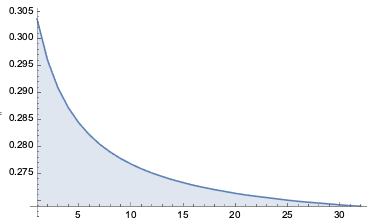}\\
\caption{Graph of $e_4(\mu)$ for $\mu\leq 36$.\label{fpro4}}
\end{center}
\end{figure}
We now show the following convergence:   
   \begin{lem}\label{hermfact1} The Fourier transform of $k_\lambda$ converges, when $\lambda\to \infty$, towards the $\Xi$-function of Riemann uniformly on closed substrips of the open  strip $\vert \Im(z)\vert<\frac 12$.
\end{lem}
\proof 
We now investigate what happens when we apply the map $\cE$ while the variable is restricted to the interval $[\lambda^{-1},\lambda]$. For $u$ in this interval, the number of integers $n$ such that $nu\leq \lambda$ is at most $\lambda/u$, thus with $$\delta(\lambda):=\max_{x\in [-\lambda,\lambda]}\vert h_\lambda(x)-h(x)\vert$$ one gets, using the definition of $\cE$ which involves $u^{1/2}$ 
$$
\vert\cE(h_\lambda)(u)-\cE(h)(u)\vert\leq u^{1/2}\delta(\lambda)\frac{\lambda}{u}
$$
 We now evaluate the Mellin transform of $k_\lambda$ on the critical strip, \ie 
 $$\cM(k_\lambda)(s)=\int_0^{\infty} u^{s-1} k_\lambda(u) du, \  \ \Re(s)\in [-\frac 12, \frac 12]$$
 We use the following estimate, where the exponent $-2$ in  $u^{-2}$ comes from two sources, and $\alpha=\Re(s)$
 $$
\vert \cM(k_\lambda)(s)-\int_{\lambda^{-1}}^\lambda k(u)u^{s-1}du\vert\leq \lambda \delta(\lambda) \int_{\lambda^{-1}}^\lambda u^\alpha u^{1/2}u^{-2}du
 $$
 One has 
 $$
 \int_{\frac{1}{\lambda }}^{\lambda } \frac{u^{\alpha +\frac{1}{2}}}{u^2} \, du=\frac{2 (\lambda ^{\frac{1}{2}-\alpha }-\lambda ^{\alpha -\frac{1}{2}})}{1-2 \alpha}
 $$
 and since $\alpha\in (-\frac 12, \frac 12)$ one has $1-2 \alpha>0$, $\frac{1}{2}-\alpha>\alpha -\frac{1}{2}$ which gives using \eqref{esti1}, \ie $\delta(\lambda)\leq c\, \lambda^{-2}$
 $$
\vert \cM(k_\lambda)(s)-\int_{\lambda^{-1}}^\lambda k(u)u^{s-1}du\vert\leq 2c \lambda ^{-1} \lambda ^{\frac{1}{2}-\alpha }(1-2 \alpha)^{-1}=
2c  \lambda ^{-\frac{1}{2}-\alpha }(1-2 \alpha)^{-1}
 $$
 Hence, since $\alpha\in (-\frac 12, \frac 12)$, one has $\frac 12 +\alpha>0$ and one  obtains, for fixed $\alpha$ 
 $$
 \vert \cM(k_\lambda)(s)-\int_{\lambda^{-1}}^\lambda k(u)u^{s-1}du\vert=O(\lambda ^{-\frac{1}{2}-\alpha })
 $$ 
 It remains to control the remainder in the Mellin transform of $k$. By the Poisson formula one has $k(u)=k(u^{-1})$ and thus it is enough to control 
 $$
 \int_\lambda^\infty k(u)u^{s-1}du
 $$
 but this tends to $0$ when $\lambda\to \infty$ due to the convergence of the integral.\endproof

\section{The missing steps}\label{missing}

There are two essential steps still missing to justify our tentative proof of the Riemann Hypothesis. 
The first is that, in order to apply Theorem~5.10 to the Weil quadratic form \( QW_\lambda \), one must prove that its smallest eigenvalue---whose existence is ensured by Theorem \ref{thmsmallest}---is simple and that its corresponding eigenvector \( \xi_\lambda \) is even. 
The second step is to establish that \( k_\lambda \) provides a sufficiently accurate approximation to (a scalar multiple of) \( \xi_\lambda \), in order to justify the convergence of the zeros of $\widehat \xi_\lambda$ towards the non-trivial zeros of \( \zeta(\frac 12 +is) \). 

\begin{figure}[H]
\begin{center}
\includegraphics[scale=0.37]{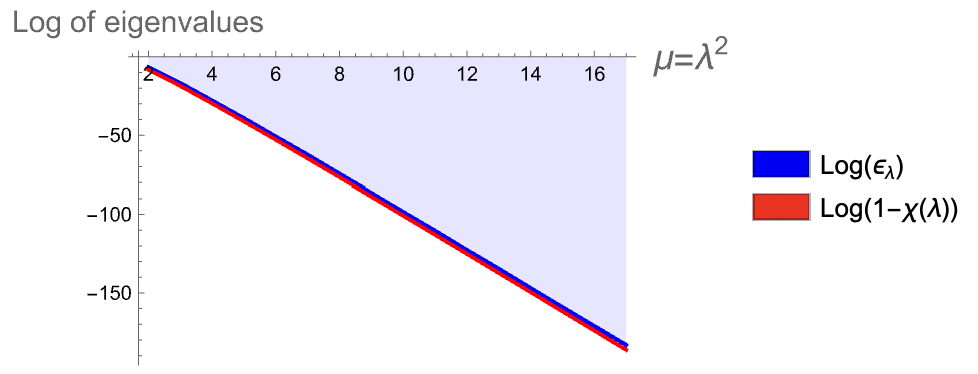}\\
\caption{Graphs of $\log(\epsilon_\lambda))$ and $\log(1-\chi(\lambda)))$ as functions of $\mu=\lambda^2$.\label{fpro1}}
\end{center}
\end{figure}

There are, however, three indications supporting the feasibility of these steps. \newline
(1) The ``simple-even'' condition holds for all values of \( \lambda \) for the prolate-wave operator \( PW_\lambda \). \newline
(2) The extremely small numbers \( \varepsilon_\lambda \) that occur as eigenvalues of the Weil quadratic form \( QW_\lambda \) also appear---see Figure \ref{fpro1}---when evaluating the discrepancy for \( h_\lambda \) to belong simultaneously to \( P_\lambda \) and \( \widehat{P}_\lambda \). \newline
(3) The numerical evidence for the proximity between \( k_\lambda \) and \( \xi_\lambda \) extends to the higher eigenfunctions of the Weil quadratic form. \newline

It remains possible that our strategy for proving convergence towards the zeros of \(  \zeta(\frac 12 +is)  \) will face significant obstacles. 
Nevertheless, it provides a strong motivation to further develop the relationship, first uncovered in~\cite{Co-zeta}, between the Weil quadratic form \( QW_\lambda \) and the prolate-wave operator \( PW_\lambda \). 
The cornerstone of this development is the trace formula established in~\cite{Co-zeta}, which relates \( P_\lambda \), \( \widehat{P}_\lambda \), and the map \( \mathcal{E} \) to \( QW_\lambda \). 
Yet, regardless of how far one can progress along this path, the present approach naturally opens the way to a deeper exploration of the unexpected relationship between two seemingly distant mathematical worlds.

\medskip
\noindent
\textbf{The world of the Weil quadratic form.} 
A key discovery of Andr\'e Weil  is the following remarkable fact: the Riemann Hypothesis is equivalent to the positivity of certain quadratic forms that involve only finitely many primes. 
This is striking, since one might expect that addressing the Riemann Hypothesis would require control over the entire infinite set of primes. 
Here, however, the problem acquires a local character, reducing to finite collections at a time. 
Moreover, in this framework, as we have seen, one can exploit our general construction of functions whose zeros lie entirely on the critical line.

\medskip
\noindent
\textbf{The world of prolate wave functions.} 
Developed by David~Slepian and collaborators, and rooted in Claude~Shannon’s work on communication theory, this theory exhibits  the miraculous relationship between the orthogonal projections that define time and frequency limitations in signal analysis and a classical second-order differential operator on the real line: the prolate wave operator, itself obtained as a confluence from the Heun equation—an object entirely familiar within Riemann’s mathematical universe. 

\medskip
\noindent
The prolate operator  plays a dual role. 
In the infrared regime, it provides a tool for approximating the minimal eigenvector of the Weil quadratic form. 
At the opposite, ultraviolet, end it also furnishes—cf.~\cite{CM}—a model of a self-adjoint operator whose spectrum reflects the high-frequency (ultraviolet) behavior of the zeros of the Riemann zeta function. 
This duality emphasizes how ideas from information theory, spectral analysis, and number theory converge within a unified framework, turning the relation between \( QW_\lambda \) and \( PW_\lambda \) into a fertile ground for further exploration.





\vspace{20pt}
\scriptsize
\noindent

\begin{tabular}{ l l l l l l l l l l l l l}
 Alain Connes &&&&&&  Caterina Consani  &&&&&&  Henri Moscovici\\ 
 Coll\`ege de France &&&&&&  Department of Mathematics  &&&&&&  Department of Mathematics\\  3 Rue d'Ulm &&&&&&  Johns Hopkins University  &&&&&& Ohio State University\\
 75005 Paris &&&&&&  Baltimore, MD 21218 &&&&&&  Columbus, OH 43210\\
 France &&&&&&  USA &&&&&&  USA\\
 \href{mailto:alain@connes.org}{alain@connes.org} &&&&&&  \href{mailto:cconsan1@jhu.edu}{cconsan1@jhu.edu}&&&&&&  \href{mailto:moscovici.1@osu.edu}{moscovici.1@osu.edu}
\end{tabular}

\end{document}